\documentclass[12pt]{article}

\UseRawInputEncoding

\usepackage{mathtools,amsthm,amssymb,color,cases}

\usepackage{stackengine}

\usepackage{mathrsfs}
\usepackage[abbrev]{amsrefs}
\usepackage{tikz}
\usepackage{enumitem}
\usepackage{setspace}

\usepackage{latexsym}
\usepackage{mathrsfs}
\usepackage{comment}
\usepackage{microtype}

\usepackage{indentfirst}

\usepackage{geometry}
\usepackage{etoolbox}
\AtEndEnvironment{proof}{\setcounter{proofpart}{0}}

 \newtheoremstyle{mystyle}% % Name
    {}%                      % Space above
    {}%                      % Space below
    {\normalfont}%           % Body font
    {}%                      % Indent amount
    {\bfseries}%             % Theorem head font
    {}%                      % Punctuation after theorem head
    { }%                     % Space after theorem head, ' ', or \newline
    {}%                      % Theorem head spec (can be left empty, meaning `normal')
  \theoremstyle{mystyle}     

\allowdisplaybreaks[1]
\theoremstyle{theorem}
\newtheorem{thm}{Theorem}[section]
\newtheorem{cor}[thm]{Corollary}
\newtheorem{lem}[thm]{Lemma}
\newtheorem*{rem*}{Remark}

\newtheorem{prop}[thm]{Proposition}

\theoremstyle{definition}
\newtheorem*{definition}{Definition}

\newtheorem{rem}[thm]{Remark}

\newtheoremstyle{part}{}{}{\normalfont}{}{\itshape}{.}{.5em}{}
\theoremstyle{part}

\newcommand\blfootnote[1]{%
  \begingroup
  \renewcommand\thefootnote{}\footnote{#1}%
  \addtocounter{footnote}{-1}%
  \endgroup
}

\numberwithin{equation}{section}
\numberwithin{thm}{section}

\newenvironment{equ*}{
    \begin{equation*}
}{
    \end{equation*}
}

\newenvironment{pf}
   {{\noindent \bf Proof.}}{\hfill \qed}

\newtheoremstyle{part}{}{}{\normalfont}{}{\itshape}{.}{.5em}{}
\theoremstyle{part}

\newcommand{\ddj}{\dot{\Delta}_j}

\newcommand{\pt}{\partial}

\newcommand{\F}{\mathcal {F}}

\newcommand{\Z}{\mathbb{Z}}

\newcommand{\R}{\mathbb R}

\newcommand{\intr}{\int_{\R^d}}

\newcommand{\al}{\alpha}

\newcommand{\x}{\xi}

\newcommand{\N}{\nabla }

\renewcommand{\div}{{\rm {div}}}

\newcommand{\supp}{\text{ supp }}

\renewcommand{\P}{\mathcal{P}}
\newcommand{\Q}{\mathcal{Q}}

\newcommand*{\dd}{\mathop{}\!\mathrm{d}}

\newcommand{\hu}{\hat{u}}
\newcommand{\hp}{\hat{\phi}}

\DeclarePairedDelimiter{\norm}{\lVert}{\rVert}
\DeclarePairedDelimiter{\bignorm}{\bigg\|}{\bigg\|}

\DeclarePairedDelimiter{\bigc}{\big{(}}{\big{)}}
\DeclarePairedDelimiter{\Bigc}{\Big{(}}{\Big{)}}
\DeclarePairedDelimiter{\Biggc}{\Bigg{(}}{\Bigg{)}}
\DeclarePairedDelimiter{\biggc}{\bigg{(}}{\bigg{)}}



\newcommand{\db}{\dot{B}}
\newcommand{\fdp}{ \frac{d}{p} }
\newcommand{\fdtwo}{ \frac{d}{2} }

\newcommand{\da}[1]{\delta a^{#1}}
\newcommand{\du}[1]{\delta u^{#1}}

\newcommand{\weakstar}{\overset{\ast}{\rightharpoonup}}

\begin{document}

\begin{center}
    {\bf \large Global existence of solutions to the compressible Navier-Stokes equations in weighted besov spaces} \\
    
    \phantom{}
    
    %Tsukasa Iwabuchi* \quad 
    D\'aith\'i \'O hAodha* \\
    
     \phantom{}
    
    Mathematical Institute, Tohoku University,
    
    980-0845
\end{center}

\blfootnote{*Email:
%*t-iwabuchi@tohoku.ac.jp,
david.declan.hughes.p6@dc.tohoku.ac.jp}

\begin{center}
\begin{minipage}{135mm}
\footnotesize
{\sc Abstract. }
We prove global existence of solutions to the Cauchy problem for the compressible Navier-Stokes equations in Euclidean spaces, given initial data with small norms in 
Besov and critical weighted Besov spaces. Global existence and a priori bounds in these spaces are obtained for dimension $d\geq 3$.
\end{minipage}
\end{center}

\section{Introduction}

In this paper, we consider the barotropic compressible Navier-Stokes system
\begin{align} \label{CNSo}
    \begin{cases}
    \pt_t \rho + \div (\rho u) = 0, & \text{ in } (0,\infty) \times\R^d, \\
    \pt_t (\rho u) + \div (\rho u \otimes u) - \div(2 \mu D(u) + \lambda \div (u) \text{Id}) + \nabla \Pi = 0, & \text{ in } (0,\infty) \times\R^d, \\
    (\rho, u) |_{t=0} = (\rho_0, u_0), & \text{ in } \R^d,
    \end{cases}
\end{align}
where $\rho : [0,\infty) \times \R^d \to [0,\infty),$ and $u: [0,\infty) \times \R^d \to \R^d$ are unknown functions, representing respectively the density and velocity of a fluid which fills all of Euclidean space $\R^d$, $d\geq 2$. 
$\Pi: [0,\infty) \times \R^d \to \R$ is the pressure in the fluid, and the barotropic assumption states that $\Pi \coloneqq P(\rho)$, for some smooth function $P(\cdot)$.
$\mu,\lambda$ are viscosity coefficients, satisfying
\begin{align} \label{coeff cndtn}
    \mu > 0, \quad 2\mu + \lambda > 0.
\end{align}
The deformation tensor $D(u)$ is defined as follows:
\begin{align*}
    D(u) \coloneqq \frac{1}{2} \Big{(} Du + Du^T \Big{)}.
\end{align*}
Where $Du$ is the matrix of derivatives of $u$, and $Du^T$ is the transpose of that matrix.
%Write what D and ^T are!

In what follows, we will show that the above problem has a unique global-in-time solution $(\rho,u)$ in critical Besov and critical weighted Besov spaces.
This result is used by the author and Iwabuchi in \cite{ohaodha2023-nonlin-est} to obtain sharp time-decay estimates for the curl-free problem.

We make some additional assumptions.
Firstly, $(\rho_0 - 1,u_0)$ is taken to be sufficiently small in the above mentioned spaces, and we define $a \coloneqq \rho - 1$.
Secondly, we assume that $\mu, \lambda$ are constant.
This allows us to rewrite the system \eqref{CNSo} as follows:
\begin{align}
\label{31}
\displaystyle
    \begin{cases}
    \displaystyle
    \pt_t a + u \cdot \N a = -(1+a) \div(u) & \text{ in } (0,\infty) \times \R^d, \\
    \displaystyle
    \pt_t u - \mathcal{A} u = - u\cdot\N u - \frac{a}{1+a} \mathcal{A} u - \frac{P'(1+a)}{1+a} \N a & \text{ in } (0,\infty) \times \R^d, \\
    \displaystyle
    (a,u)\Big{|}_{t=0} = (a_0,u_0) & \text{ in }\R^d,
    \end{cases}
\end{align}
where $\mathcal{A}$ 
is defined as follows:
\[ \mathcal{A}u \coloneqq \mu \Delta u + (\lambda + \mu) \nabla \div(u). \]

We introduce some previous results. Matsumura-Nishida showed in~\cites{matsumura-nishida1979, matsumura-nishida1980} that \eqref{CNSo} has global solutions in the $d=3$ case when equipped with data $(\rho_0, u_0)$ that is a small perturbation in $L^1 \cap H^3$ of $(\bar{\rho}, 0)$ for any positive constant $\bar{\rho}$, and proved the following decay result
\begin{align*}
    \bignorm{ 
    \begin{bmatrix}
        \rho(t) - \bar{\rho} \\
        u(t)
    \end{bmatrix}  
    }_{2}
    \leq C ({1+t})^{-3/4}.
\end{align*}
This is the decay rate of the solution to the heat equation with initial data in $L^1$.
Ponce in \cite{ponce1985} then extended these results. In particular, for $p \in [2,\infty]$, $k\in\{0,1,2\}$, and dimension $d=2,3$,
\begin{align*}
    \bignorm{ 
    \nabla^k \begin{bmatrix}
        \rho(t) - \bar{\rho} \\
        u(t)
    \end{bmatrix}  
    }_{p}
    \leq C ({1+t})^{-\frac{d}{2}(1-\frac{1}{p}) - \frac{k}{2} }.
\end{align*}

Before introducing the next results, we introduce the concepts of scaling invariance and critical spaces.
We observe that, if $(\rho,u)$ is a solution to the problem \eqref{CNSo} with initial data $(\rho_0,u_0)$, then the rescaled pair
\begin{align*}
    ( \rho(\ell^2t,\ell x) , \ell u(\ell^2t, \ell x) )
\end{align*}
is also a solution for all $\ell>0$, with the rescaled initial data $(\rho_0(\ell x), \ell u_0(\ell x))$ if we replace the pressure function $P$ with $\ell^2 P$.
We thus say that \eqref{CNSo} is \textit{invariant} under the above scaling transformation.
We then say that a function space $F \subseteq \mathcal{S}' \times (\mathcal{S}')^d$ is \textit{critical} if its associated norm is invariant under the transformation 
\[
(\rho,u) \mapsto ( \rho(\ell \cdot) , \ell u(\ell \cdot) ).
\]
We then observe that the space $\db^\fdp_{p,1} \times \biggc{\db^{\fdp-1}_{p,1}}^d$ is critical in this sense.

Moving onto results in critical spaces, Danchin first solved the global existence problem for \eqref{CNSo} with small initial data $(\rho - \bar{\rho} , u_0)$ in the critical Besov space $\dot{B}^{ \frac{d}{2} }_{2,1} \times \dot{B}^{ \frac{d}{2} - 1}_{2,1}$ in~\cite{Danchin2000}. 
Danchin also analyses this problem in~\cite{danchin2016}, with the high-frequency part of the norm generalised to $p$ close to $2$.
In both, an additional assumption that the low-frequency part of $\rho - \bar{\rho}$ be small also in $\dot{B}^{ \frac{d}{2} - 1}_{2,1}$ was required for global existence, while Danchin demonstrates in \cite{danchin2016} that, for $p \in [2,2d)$, smallness of $a_0$ in the critical space $\dot{B}^{  \fdp }_{p,1}$ and existence of $u_0$ in $\dot{B}^{ \fdp  - 1}_{p,1}$ alone are sufficient for local-in-time existence.
Time-decay estimates in Besov norms for $p$ close to $2$ were proven by Danchin-Xu in~\cite{Danchin-Xu2017}.
In this paper, we will often speak of `\textit{subcritical}' spaces and their norms, which we define simply as those Besov spaces with regularity exponents higher than those of the critical space.

Our main result 
uses the results in the critical Besov framework in order to obtain global existence of solutions in weighted Besov spaces. 
Before we give our main result, we give the definitions of Besov spaces, weighted Besov spaces, and high-frequency and low-frequency norms.

\begin{definition} \label{besov spaces}
Let $\{ \hat{\phi}_j \}_{ j \in\Z}$ be a set of non-negative measurable functions such that 
\begin{enumerate}
    \item $\displaystyle \sum_{ j \in\Z} \hat{\phi}_j (\x) = 1, \text{ for all } \x \in \R^3 \backslash \{0\}$,
    \item $\hat{\phi}_j (\x) = \hat{\phi}_0(2^{-j}\x)$,
    \item $\supp \hat{\phi}_j (\x) \subseteq \{ \x \in \R^3 \ | \ 2^{j-1} \leq |\x| \leq 2^{j+1} \}$.
\end{enumerate}
For a tempered distribution $f$, we write
\[
\dot{\Delta}_j f \coloneqq \F^{-1} [\hat{\phi}_j \hat{f}].
\]
This gives us the \textit{Littlewood-Paley decomposition} of $f$:
\begin{align*}
    f = \sum_{j\in\Z} \ddj f.
\end{align*}
This equality only holds modulo functions whose Fourier transforms are supported at $0$, i.e. polynomials.
To ensure equality in the sense of distributions, we
next let $\dot{S}_j$ denote the sum of dyadic blocks up to $j$.
That is, for $j\in\Z$,
\begin{align*}
    \dot{S}_j f \coloneqq \sum_{j'\leq j} \dot{\Delta}_{j'} f.
\end{align*}
Then we consider the subset $\mathcal{S}'_h$ of tempered distributions $f$ such that
\begin{align*}
    \lim_{j\to -\infty} \norm{ \dot{S}_j f }_{L^\infty} = 0.
\end{align*}
The Besov norm is then defined as follows: for $1 \leq p,q \leq \infty$, and $s \in \R$, we define 
\[
\norm{f}_{\dot{B}^{s}_{p,q}} \coloneqq \Bigc{ \sum_{j \in \Z} 2^{sqj} \norm{\dot{\Delta}_j f}^q_{p} }^{ \frac{1}{q} }.
\]
The set $\dot{B}^{s}_{p,q}$ is defined as the set of functions, $f \in \mathcal{S}'_h$, whose Besov norm is finite. 
Throughout this paper, we will refer to the parameter $s$ as the `\textit{regularity exponent}' and $p$ as the `\textit{Lebesgue exponent}.'

We then define a weighted Besov space as the set %$\dot{B}^{s}_{p,q}(x_k)$ 
of functions $f \in \mathcal{S}'_h$ such that the Besov norm of $f$ multiplied by $x_k$ is finite for all $k \in \{ 1, 2, ... , d \}$. That is,
\[
\norm{x_k f}_{\dot{B}^{s}_{p,q}} < \infty \text{ for all } k \in \{ 1, 2, ... , d \}. 
\]
We then call $\norm{x_k f}_{\dot{B}^{s}_{p,q}}$ the weighted Besov norm of $f$.

We also regularly use the following notation for so-called high-frequency and low-frequency norms:
\[
\norm{f}^h_{\dot{B}^{s}_{p,q}} \coloneqq \Bigc{ \sum_{j \geq j_0} 2^{sqj} \norm{\dot{\Delta}_j f}^q_{p} }^{ \frac{1}{q}}, \quad
\norm{f}^l_{\dot{B}^{s}_{p,q}} \coloneqq \Bigc{ \sum_{j \leq j_0} 2^{sqj} \norm{\dot{\Delta}_j f}^q_{p} }^{ \frac{1}{q}},
\]
where $j_0 \in \Z$ is called the frequency cut-off constant.
We also define the high-frequency and low-frequency parts of a function $f$:
\[
f^h \coloneqq \sum_{j \geq j_0} \dot{\Delta}_j f, 
\quad
f^l \coloneqq \sum_{j \leq j_0} \dot{\Delta}_j f.
\]
\end{definition}

Finally, we introduce our main result. Following Danchin, we take the low-frequency part of initial data $(a_0,u_0)$ to be small in $\dot{B}^{ \frac{d}{2} - 1}_{2,1}$. That is, the low-frequency part of $a_0$ requires some additional regularity, being in the space `one derivative less than' the critical space.
The high-frequency part of $(a_0,u_0)$ is then taken small in the subcritical space $\dot{B}^{\frac{d}{2}+1}_{2,1} \times \dot{B}^{ \frac{d}{2} }_{2,1}$, where the regularity exponents of the spaces containing $a_0$ and $u_0$ are one higher than their respective critical spaces. We shall see in section \ref{section local} that this is necessary in order to obtain local existence in weighted spaces.

\newpage

\begin{definition}
    We introduce the solution space $S$ to which our solution will belong as the set of all pairs of functions $(a,u)$, where $a : [0,\infty) \times \R^d \to [0,\infty)$ is a scalar function and $u: [0,\infty) \times \R^d \to \R^d$ is a $d$-vector function, satisfying the following for all $k \in \{ 1, 2, ..., d \}$:
\begin{align*}
    & (a,u)^l \in \widetilde{C} (\R_{>0}; \db^{\frac{d}{2} - 1}_{2,1}) \cap L^1(\R_{>0}; \db^{\frac{d}{2} + 1}_{2,1}), 
    \quad
    a^h \in \widetilde{C} (\R_{>0}; \db^{\frac{d}{2} + 1}_{2,1}) \cap L^1(\R_{>0}; \db^{\frac{d}{2} + 1}_{2,1}),
    \\
    & u^h \in \widetilde{C}(\R_{>0}; \db^{\frac{d}{2}}_{2,1})
    \cap L^1(\R_{>0}; \db^{\frac{d}{2} + 2}_{2,1}),
    \\
    & (x_k a, x_k u)^l \in \widetilde{C} (\R_{>0}; \db^{\frac{d}{2}}_{2,1}) 
    \cap L^1_t ( \R_{>0} ; \db^{\frac{d}{2} + 2}_{2,1} ), 
    \quad 
    (x_k a)^h \in \widetilde{C} (\R_{>0}; \db^{\frac{d}{2} +1 }_{2,1})
    \cap L^1_t ( \R_{>0} ; \db^{\frac{d}{2} + 1}_{2,1} ),
    \\
    & 
    % (x_k u)^l \in \widetilde{C} (\R_{>0}; \db^{\frac{d}{2}}_{2,1})
    % \cap L^1_t ( \R_{>0} ; \db^{\frac{d}{2} + 2}_{2,1} ),
    % \quad
    (x_k u)^h \in \widetilde{C} (\R_{>0}; \db^{\frac{d}{2}}_{2,1})
    \cap L^1_t ( \R_{>0} ; \db^{\frac{d}{2} + 2}_{2,1} ).
\end{align*}
$S$ is then equipped with the obvious norm corresponding to the strong topologies for the above spaces.
\end{definition}

Here, $\widetilde{C} (\R_{>0}; \db^{s}_{p,1}) \coloneqq 
{C} (\R_{>0}; \db^{s}_{p,1})
\cap
\widetilde{L}^\infty(\R_{>0}; \db^{s}_{p,1}),
$ for $s\in\R, \ p \in[1,\infty].$ 
The norm of $\widetilde{L}^\infty( 0,T ; \db^{s}_{p,1})$ for $T>0$ is defined by taking the $L^\infty$-norm over the time interval \textit{before} summing over $j$ for the Besov norm. That is, for all $f \in \widetilde{L}^\infty( 0,T ; \db^{s}_{p,1})$,
\[
\norm{f}_{ \widetilde{L}^\infty( 0,T ; \db^{s}_{p,1}) }
\coloneqq
\sum_{j\in\Z} 2^{sj} \sup_{ t\in(0,T) } \norm{ \ddj f(t) }_{L^p} .
\]
We abbreviate the notation for norms by writing 
\[
\norm{f}_{ \widetilde{L}^\infty_T \db^{s}_{p,1} }
\coloneqq
\norm{f}_{ \widetilde{L}^\infty( 0,T ; \db^{s}_{p,1}) },
\]
and similarly abbreviate other norms over time and space.

% \begin{thm} 
%     Let $d\geq 3,$ $k \in \{ 1, 2, ..., d \},$ and $p \in [2, \text{min} \{ d, \, 2d/(d-2) \} )$. 
%     Assume $P'(1) > 0.$
%     Then there exists 
%     a frequency cut-off constant $j_0 \in \Z$ and
%     %%%%%%%%%%
%     %Move remark abt the freq cutoff being large enough to the proof?
%     %%%%%%%%%%
%     a small constant $c = c(p,d,\mu,P) \in \R$ such that, if $(a_0,u_0)$ 
%     satisfy 
%     \begin{align}
%     \notag
%         S_{p,0} \coloneqq &  
%         \norm{(a_0,u_0)}^l_{\dot{B}^{\frac{d}{2} - 1}_{2,1}} 
%         + \norm{a_0}^h_{\db^{\frac{d}{p} + 1}_{p,1}} 
%         + \norm{u_0}^h_{\db^{\frac{d}{p}}_{p,1}}
%         \\
%         & + \norm{ ( x_k a_0, x_k u_0 ) }^l_{\db^{\frac{d}{2}}_{2,1} }
%         +
%         \norm{x_k a_0}^h_{ \db^{\frac{d}{p} + 1}_{p,1}}
%         + \norm{x_k u_0}^h_{\db^{\frac{d}{p}}_{p,1}}
%         +
%         \norm{(a_0,u_0)}_{ \db^{-s_0}_{2,\infty} }
%         \leq c,
%     \end{align}
%     where $s_0 \coloneqq d (\frac{2}{p} - \frac{1}{2})$,
%     then \eqref{31} has a unique global-in-time solution $(a,u)$ in the space $S_p$ defined above. Also, there exists a constant $C = C(p,d,\mu, P, j_0)$ such that
%     \begin{align}
%         \norm{(a,u)}_{S_p} \leq C S_{p,0}.
%     \end{align}
% \end{thm}

\begin{thm} 
    Let $d\geq 3$. 
    %$k \in \{ 1, 2, ..., d \}$. 
    Assume $P'(1) > 0.$
    Then there exists 
    a frequency cut-off constant $j_0 \in \Z$ and
    a small constant $c = c(d,\mu,P) \in \R$ such that, if $(a_0,u_0)$ 
    satisfy 
    \begin{align*}
    \notag
        S_{0} \coloneqq &  
        \norm{(a_0,u_0)}^l_{\dot{B}^{\frac{d}{2} - 1}_{2,1}} 
        + \norm{a_0}^h_{\db^{\frac{d}{2} + 1}_{2,1}} 
        + \norm{u_0}^h_{\db^{\frac{d}{2}}_{2,1}}
        \\
        & + \sum_{k=1}^d
         \Bigc{\norm{ ( x_k a_0, x_k u_0 ) }^l_{\db^{\frac{d}{2}}_{2,1} }
        +
        \norm{x_k a_0}^h_{ \db^{\frac{d}{2} + 1}_{2,1}}
        + \norm{x_k u_0}^h_{\db^{\frac{d}{2}}_{2,1}} }
        +
        \norm{(a_0,u_0)}_{ \db^{-\frac{d}{2}}_{2,\infty} }
        \leq c,
    \end{align*}
    then \eqref{31} has a unique global-in-time solution $(a,u)$ in the space $S$ defined above. Also, there exists a constant $C = C(d,\mu, P, j_0)$ such that
    \begin{align*}
        \norm{(a,u)}_{S} \leq C S_{0}.
    \end{align*}
\end{thm}

\begin{rem}
    The smallness of the 
    %$\db^{-s_0}_{2,\infty}$-norm 
    $\db^{-\frac{d}{2}}_{2,\infty}$-norm 
    is required in order to take advantage of a low-frequency estimate due to Danchin and Xu in \cite{Danchin-Xu2017}, which we use to extend the local weighted solutions $(A_k,U_k)$ to global solutions in Section \ref{section global}. 
\end{rem}

We briefly discuss the proof.
We show unique local existence on a time interval $[0,T]$, with $T>0,$ of a solution $(a,u)$ to the above simplified system in Besov spaces using an argument by convergence from an approximate system. 
In order to obtain local existence in weighted Besov spaces, we multiply the entire system \eqref{31} by $x_k,$ where $k\in\{1,2,...,d\}$; and then, setting $A_k \coloneqq x_k a,$ $U_k \coloneqq x_k u,$ we apply the same steps to the new system
\begin{align}
\label{31x}
\displaystyle
    \begin{cases}
    \displaystyle
    \pt_t A_k + u \cdot \N A_k = -( x_k + A_k ) \div(u) + au_k & \text{ in } (0,\infty) \times \R^d, \\
    \displaystyle
    \pt_t U_k - \mathcal{A} U_k =  
    - U_k \cdot\N u 
    - \frac{ a }{1+a} \mathcal{A} U_k 
    - \frac{P'(1+a)}{1+a} 
    \biggc{
    \N A_k - a e_k
    }
    \\
    \displaystyle \quad \quad \quad \quad \quad \quad \quad  
    +\frac{1}{1+a}
    \biggc{
    2\mu \pt_k u
    +(\lambda +\mu) \bigc{ \div(u) e_k + \N(u_k) }
    }
    & \text{ in } (0,\infty) \times \R^d, \\
    \displaystyle
    (a,u)\Big{|}_{t=0} = (a_0,u_0) & \text{ in }\R^d,
    \end{cases}
\end{align}
and obtain a unique local solution $(A_k,U_k)$.
Here, $e_k$ denotes the unit vector along the $x_k$-axis.
The local solutions are then extended past $T$ by proving a priori estimates for both systems.
We will see on page \pageref{reason hi reg Uk} in the proof for local existence that the term $\frac{P'(1+a)}{1+a} a e_k$
forces us to add extra regularity on the high frequencies of 
$(a_0,u_0)$, compared to the local existence proof found in \cite{danchin2016}.

%%%%=====================================================================================================================================================================================================================================================================================================

%PRELIMINARIES

%%%%=====================================================================================================================================================================================================================================================================================================

%%%Add a list of notation?

\section{Preliminaries}

% \subsection{$L^p$ Spaces}

% \begin{definition} ($L^p$ Space)

% Let $(\Omega, M, \mu)$ be a measure space. The set $L^p(\Omega)$, with $1\leq p<\infty$, is defined as the set of all functions, $f:\Omega \to \R$, such that
% \[
% \norm{f}_{p} := \Bigc{\int_\Omega |f(x)|^p \dd \mu }^{1/p} <\infty.
% \]

% For $p=\infty$, the set is defined as the set of all functions whose essential supremum is finite.
% \[
% \norm{f}_{\infty} := \inf \{ C\geq0 \ | \ |f(x)| \leq C, \text{ for almost all } x \in \Omega \}.
% \]

% Functions that agree almost everywhere (i.e. everywhere in $\Omega$ except on a subset with $0$ measure) are considered a single element of $L^p(\Omega)$. 

% \end{definition}

% \noindent \textbf{Notation.} For $1\leq p \leq \infty$, we denote by $p'$ the conjugate exponent. That is,
% \[
% \frac{1}{p} + \frac{1}{p'} = 1,
% \]
% with the convention that $\frac{1}{\infty} := 0$ in this context.

% \begin{prop}(H\"older's Inequality)

% Let $f \in L^p(\Omega)$ and $g \in L^{p'}(\Omega)$, where $1 \leq p \leq \infty$. Then $fg \in L^1$, and
% \[
% \int_\Omega |fg| \dd\mu \leq \norm{f}_{p} \norm{g}_{p'}.
% \]

% \end{prop}

% \begin{prop}(Young's Convolution Inequality)

% Let $f \in L^p$ and $g \in L^q$, where $1 \leq p,q \leq r \leq \infty$, such that 
% \[
% \frac{1}{p} + \frac{1}{q} = \frac{1}{r} + 1.
% \]
% Then
% \[
% \norm{f\ast g}r \leq \norm{f}_{p}\norm{g}_{q}.
% \]
% Here, $(f \ast g)$ denotes the convolution of $f$ and $g$. That is,
% \[
% (f \ast g)(x) := \int_{\R^2} f(x-y) g(y) \dd y = \int_{\R^2} f(y) g(x-y) \dd y.
% \]

% \end{prop}

\begin{definition}(The Fourier Transform)
For a function, $f$, we define the Fourier transform of $f$ as follows:
\begin{align*}
\F[f](\x) \coloneqq \hat{f}(\x) \coloneqq \frac{1}{(2\pi)^{d/2}}\intr e^{-i x\cdot\x} f(x) \dd x.
\end{align*}
The inverse Fourier transform is then defined as 
\begin{align*}
\F^{-1}[\hat{f}](x) \coloneqq \frac{1}{(2\pi)^{d/2}} \intr e^{i x \cdot \xi} \hat{f}(\x) \dd \x.
\end{align*}
For the purpose of calculating inequalities, we will frequently omit the factor of $1/(2\pi)^{d/2}$.
\end{definition}

% \begin{prop}(Plancherel Theorem)

% Let $f \in L^1 \cap L^\infty.$ Then the $L^2$ norm of $f$ is invariant under the Fourier transform. That is,
% \[
% \norm{f}_2 = \norm{\hat{f}}_2.
% \]

% \end{prop}

\begin{definition}(Orthogonal Projections on the divergence-free and curl-free fields)
% We define the projection mappings $\P, \Q$ using the Fourier transform and the Kronecker delta, defined as follows: 
% \[ \delta_{i,j} \coloneqq  
% \begin{cases}
% 0, & i \neq j, \\
% 1, & i = j.
% \end{cases} 
% \]
The projection mapping $\P$ is a matrix with each component defined as follows for $i,j \in \{1, 2, ... , d\}$:
\begin{align*}
(\P)_{i,j} \coloneqq \delta_{i,j} + (-\Delta)^{-1} \pt_i \pt_j.
\end{align*}
We then define $\Q \coloneqq 1 - \P$.
For $f \in (\dot{B}^s_{p,q}(\R^d))^d$, with $s \in \R$, and $p, q \in [1,\infty]$, we may write
\begin{align*}
\P f \coloneqq (1 + (-\Delta)^{-1} \nabla \div) f.
\end{align*}
%where we use the Fourier transform to define the operator $(-\Delta)^{-1}f \coloneqq \F^{-1}[|\x|^{-2} \hat{f}]$.
\end{definition}

%\subsection{Besov Spaces}

% \begin{definition}

% We use the Littlewood-Paley decomposition of unity to define homogeneous Besov spaces. Let $\{ \hat{\phi}_k \}_{k\in\Z}$ be a set of non-negative measurable functions such that 
% \begin{enumerate}
%     \item $\displaystyle \sum_{k\in\Z} \hat{\phi}_k (\x) = 1, \text{ for all } \x \in \R^2 \backslash \{0\}$,
%     \item $\hat{\phi}_k(\x) = \hat{\phi}_0(2^{-k}\x)$,
%     \item $\supp \hat{\phi}_k (\x) \subseteq \{ \x \in \R^2 \ | \ 2^{k-1} \leq |\x| \leq 2^{k+1} \}$.
% \end{enumerate}
% For $f \in \mathcal{S'}/\mathcal{P}$, we write 
% \[
% \dot{\Delta}_k f \coloneqq \F^{-1} [\hat{\phi}_k \hat{f}].
% \]
% The Besov norm is then defined as follows: for $1 \leq p,q \leq \infty$, and $s \in \R$, we define 
% \[
% \norm{f}_{\dot{B}^{s}_{p,q}} := \Bigc{ \sum_{k \in \Z} 2^{sqk} \norm{\dot{\Delta}_k f}^q_{p} }^{1/q}.
% \]
% The set $\dot{B}^{s}_{p,q}$ is defined as the set of functions, $f \in \mathcal{S'}/\mathcal{P}$, whose Besov norm is finite.
% Finally, we also regularly use the following notation for so-called high and low-frequency norms:
% \[
% \norm{f}^h_{\dot{B}^{s}_{p,q}} := \Bigc{ \sum_{k \geq 3} 2^{sqk} \norm{\dot{\Delta}_k f}^q_{p} }^{1/q}, \quad
% \norm{f}^l_{\dot{B}^{s}_{p,q}} := \Bigc{ \sum_{k \leq 2} 2^{sqk} \norm{\dot{\Delta}_k f}^q_{p} }^{1/q}.
% \]
% \end{definition}

We next write some key properties of Besov spaces, whose proofs can be found in \cite{danchinbook}.

\begin{prop}
    Let $p\in[1,\infty]$.
    Then we have the following continuous embeddings:
    \begin{align*}
        \db^0_{p,1} \hookrightarrow L^p \hookrightarrow \db^0_{p,\infty}.
    \end{align*}
\end{prop}

\begin{prop}
    Let $s\in\R,$ $1\leq p_1 \leq p_2 \leq \infty$, and $1\leq r_1 \leq r_2 \leq \infty$.
    Then
    \begin{align*}
        \db^s_{p_1,r_1} \hookrightarrow \db^{s - d ( \frac{1}{p_1} - \frac{1}{p_2} )}_{p_2,r_2}.
    \end{align*}
\end{prop}

\begin{prop}
    Let $1\leq p \leq q \leq \infty$.
    Then
    \begin{align*}
        \db^{ \frac{d}{p} - \frac{d}{q} }_{p,1} \hookrightarrow L^q.
    \end{align*}
    Also, if $p < \infty$, then $\db^{ \frac{d}{p} }_{p,1}$ is continuously embedded the space $C_0$ of bounded continuous functions vanishing at infinity.
\end{prop}

For the next proposition, we introduce the notation $F(D)u \coloneqq \F^{-1} [ F \hu ]$.
\begin{prop}\label{fm est for besov} {\rm (}Fourier Multiplier Estimate{\rm )}
    Let $F$ be a smooth homogeneous function of degree $m$ on $\R^d \backslash \{0\}$ such that $F(D)$ maps $\mathcal{S}'_h$ to itself.
    Then
    \begin{align*}
        F(D) : \db^s_{p,r} \to \db^{s-m}_{p,r}.
    \end{align*}
    In particular, the gradient operator maps $\db^s_{p,r}$ to $\db^{s-1}_{p,r}$.
\end{prop}

\begin{prop} {\rm (}Composition Estimate{\rm )}
    Let $F : \R \to \R$ be smooth with $F(0) = 0$.
    Let $s>0$ and $1\leq p, r \leq \infty$.
    Then $F(u) \in \db^s_{p,r} \cap L^\infty$ for $u \in \db^s_{p,r} \cap L^\infty$, and there exists a constant $C = C( \norm{u}_{L^\infty}, F', s,p,d )>0$ such that
    \begin{align*}
        \norm{F(u)}_{ \db^s_{p,r} }
        \leq C
        \norm{ u }_{ \db^s_{p,r} }.
    \end{align*}
\end{prop}

We next introduce the Bony decomposition, which is used to study products of two tempered distributions in Besov spaces.
Let $u, v$ be two tempered distributions in
$\mathcal{S}'_h$. We have
\begin{align*}
    u= \sum_{j\in\Z} \ddj u, \text{ and } v=\sum_{j'\in\Z} \dot{\Delta}_{j'} v,
\end{align*}
and so, formally, the product may be written
\begin{align*}
    uv = \sum_{j\in\Z} \sum_{j' \in\Z} (\ddj u) (\dot{\Delta}_{j'} v).
\end{align*}
We shall split the above sum into three parts, one where $\ddj u$ is supported on large frequencies (that is, $\supp(\hp_j \hu$) is large) relative to $\dot{\Delta}_{j'} v$, vice-versa, and one final part where the frequencies are roughly the same.

\begin{definition}(Bony Decomposition)
    The \textit{homogeneous paraproduct} of $v$ by $u$ is defined:
    \begin{align*}
        T_u v \coloneqq \sum_{j\in\Z} \dot{S}_{j-4} \ddj v
    \end{align*}
    The \textit{homogeneous remainder} of $u$ and $v$ is defined:
    \begin{align*}
        R(u,v) \coloneqq \sum_{ |j-j'| \leq 3 } \ddj u \dot{\Delta}_{j'} v.
    \end{align*}
    Then, formally, we may decompose the product $uv$ as follows:
    \begin{align*}
        uv = T_u v + R(u,v) + T_v u.
    \end{align*}
\end{definition}

This decomposition allows us to bound $uv$ in Besov spaces by applying the following inequalities.
\begin{prop}
    Let $1\leq p, p_1, p_2 \leq \infty$ with $1/p = 1/p_1 + 1/p_2$. 
    Let $s\in\R$, $r\in [1,\infty]$, and $t<0$.
    Then there exists a constant $C>0$ such that
    \begin{align*}
        & \norm{T_u v}_{ \db^s_{p,r} } \leq C \norm{u}_{L^{p_1}} \norm{v}_{ \db^s_{p_2,r} },
        \\
        & \norm{ T_u v }_{ \db^{s+t}_{p,r} } 
        \leq
        C \norm{ u }_{ \db^t_{p_1,\infty} }
        \norm{v}_{ \db^s_{p_2,r} }.
    \end{align*}
    Also, if $s_1 + s_2 > 0$ and $1/r = 1/r_1 + 1/r_2$, then
    \begin{align*}
        \norm{ R(u,v) }_{ \db^{s_1 + s_2}_{p,r} }
        \leq C
        \norm{ u }_{ \db^{s_1}_{p_1,r_1} }
        \norm{ v }_{ \db^{s_2}_{p_2,r_2} }.
    \end{align*}
\end{prop}
Combining the above inequalities leads to the following important corollary, which we use regularly.
\begin{cor}
    Let $u,v \in L^\infty \cap \db^s_{p,r}$, with $s>0$ and $1\leq p, r \leq \infty$.
    Then there exists a constant $C = C(d,p,s)>0$ such that
    \begin{align*}
        \norm{uv}_{ \db^s_{p,r} }
        \leq
        C \biggc{
        \norm{ u }_{L^\infty} \norm{ v }_{ \db^s_{p,r} }
        +
        \norm{ v }_{L^\infty} \norm{ u }_{ \db^s_{p,r} }
        }.
    \end{align*}
\end{cor}

%%%%%%%%%%%%%%%%%%%%%%%%%%%%%%%%%%%
%%%%%%%%%%%%%%%%%%%%%%%%%%%%%%%%%%%
%The below Lemma is useless. We only need it in Danchin (2016) in the s=d/p+1 case, which can just be checked directly easily! Can also verify the s=d/p case.
%Danchin (2016) MISTAKENLY says he uses the s=d/p-1 case, when he actually means s=d/p+1.
%I acc have no idea what range of s this ineq holds for, but it holds for the one case we acc need so who cares.
%Also Danchin-He 2016 does NOT prove this lemma, as D claims in 2016. Iosa 
%%%%%%%%%%%%%%%%%%%%%%%%%%%%%%%%%%%
%%%%%%%%%%%%%%%%%%%%%%%%%%%%%%%%%%%

% Finally before moving on from the Bony decomposition, we give the following lemma, which is important to the proof of global existence for our main theorem. The proof can be found in \cite{danchin-he2016}.
% \begin{lem}
%     Let $d\geq 2$ and $1\leq p\leq \text{min} \{ 4, 2d/(d-2) \}$. Then
%     \begin{align}
%         \norm{ T_u v }_{ \db^{s - 1 + \frac{d}{2} - \frac{d}{p} }_{2,1} }
%         \leq C
%         \norm{ u }_{ \db^{ \frac{d}{p} - 1 }_{p,1} }
%         \norm{ v }_{ \db^s_{p,1} }.
%     \end{align}
% \end{lem}

The key to our analysis is to consider \eqref{31} and \eqref{31x} as couplings of a transport equation for the density $a, \ A_k$ and a heat equation for the velocity $u,\ U_k$.
We then may apply inequalities to the density and velocity in Besov norms for those respective equations.
We write those inequalities now, the proofs of which can be found in \cite{danchin2016}.
More detailed treatment of these equations and the compressible Navier-Stokes equations in Besov spaces can be found in \cite{danchinbook}.

%%%%%%%%%%%%%
%Rename to propositions?
%%%%%%%%%%%%%
\begin{prop} {\rm (\cite{danchin2016})} {\rm(}Besov norm estimates for the transport equation{\rm)} \label{thm2}
Let $p \in [1,\infty]$, and $s\in\R$ such that
\[
-\text{min} \{ d/p, \, d/p'\}  < 
s 
\leq 1 + d/p .
\]
    Let $a$ solve the following Cauchy problem:
    \begin{align*}
        \begin{cases}
            \pt_t a + v \cdot \N a + \lambda a = f,
            & \text{in } \R \times \R^d,
            \\
            a|_{t=0} = a_0, & \text{in } \R^d,
        \end{cases}
    \end{align*}
    where $a_0 \in \db^s_{p,1}$, $\lambda \geq 0$, $f\in L^1 (\R_{>0} ; \db^s_{p,1} )$, and $v \in L^1 (\R_{>0} ; \db^{\frac{d}{p} + 1}_{p,1} )$.
    Then, for all $t>0$,
    $a$ satisfies the following inequality:
    \begin{align*}
        \norm{a}_{ \widetilde{L}^\infty_t \db^s_{p,1} }
        +
        \lambda \norm{a}_{ L^1_t \db^s_{p,1} }
        \leq e^{C V(t) }
        \biggc{
        \norm{a_0}_{ \db^s_{p,1} }
        +
        \norm{f}_{ L^1_t \db^s_{p,1} },
        }
    \end{align*}
    with $V(t) \coloneqq \norm{\N v }_{ L^1_t \db^{\frac{d}{p}}_{p,1} }$.
\end{prop}

\begin{prop} {\rm (}\cite{danchin2016}{\rm )} (Besov norm estimates for the heat equation) \label{thm1}
    Let $p\in [1,\infty]$ and $s\in\R$.
    Let $u$ solve the following Cauchy problem:
    \begin{align*}
        \begin{cases}
            \pt_t u - \Delta u = f, & \text{in } \R_{>0} \times \R^d,
            \\
            u|_{t=0} = u_0, & \text{in } \R^d,
        \end{cases}
    \end{align*}
    where $u_0 \in \db^s_{p,1}$ and $f \in L^1 ( \R_{>0} ; \db^s_{p,1} )$. 
    Then $u$ satisfies the following inequality for all $t>0$:
    \begin{align*}
        \norm{ u }_{ \widetilde{L}^\infty_t \db^{s}_{p,1} }
        +
        \norm{ u }_{ L^1_t \db^{s+2}_{p,1} }
        \leq 
        C \biggc{
        \norm{u_0}_{\db^s_{p,1}}
        +
        \norm{f}_{ L^1_t \db^{s}_{p,1} }
        }.
    \end{align*}
\end{prop}

We treat $u$ in \eqref{31} as satisfying the Lam\'{e} equation
\begin{align}\label{Lame}
    \pt_t u - \mathcal{A} u = g,
\end{align}
where we write $g$ for the right-hand side of the momentum equation.
Then, splitting $u$ into its divergence-free and curl-free parts using the projection operators $\P$ and $\Q$, we see that each part satisfies a heat equation. 
That is,
\begin{align*}
    \begin{cases}
        \pt_t \P u - \mu \Delta \P u = \P g,
        \\
        \pt_t \Q u - \nu \Delta \Q u = \Q g.
    \end{cases}
\end{align*}
Applying Proposition \ref{thm1} to both equations, and combining the two inequalities (noting that $P$ and $Q$ are both $0$-order Fourier multipliers and applying Proposition \ref{fm est for besov}) yields
\begin{align} \label{ineq for u}
    \norm{ u }_{ \widetilde{L}^\infty_t \db^{s}_{p,1} }
        +
        \text{min} \{ \mu, \nu \}
        \norm{ u }_{ L^1_t \db^{s+2}_{p,1} }
        \leq 
        C \biggc{
        \norm{u_0}_{\db^s_{p,1}}
        +
        \norm{g}_{ L^1_t \db^{s}_{p,1} }
        },
\end{align}
and similar for $U_k$.

%%%%=====================================================================================================================================================================================================================================================================================================

%PROOF OF LOCAL EX

%%%%=====================================================================================================================================================================================================================================================================================================

\section{Local Existence and Uniqueness}
\label{section local}

In this section, we prove that \eqref{31} and \eqref{31x} have unique local solutions on a time interval $[0,T]$ for some $T>0$.
The proof for \eqref{31x} relies on that of \eqref{31}, and so we roll the existence results for both systems into one proposition.

Due to the regularity exponents in inequality \eqref{ineq for u} and the nonlinear terms of \eqref{31}, it is natural when starting with $u_0 \in \db^{\frac{d}{p} - 1}_{p,1}$ to expect $u$ to be in $E^{\frac{d}{p} - 1}_p(T)$, where $E^{ s }_p(T)$ is defined below for all $p\in [1,\infty]$, $s\in\R$, and $T>0$:
\begin{align*}
    E^{ s }_p(T)
    \coloneqq
    \{
    u \in C( [0,T] ; \db^{s}_{p,1} ) \ | \ \pt_t u, \, \N^2 u \in L^1 (0,T ; \db^{s}_{p,1}) 
    \}.
\end{align*}
We shall see, however, that additional regularity for $a_0, u_0$ in the Besov spaces with exponents one higher than their respective critical spaces is also necessary to obtain local existence for \eqref{31x}.
This is because the critical spaces for $x_k a, x_k u$ have regularity exponents one high than those for $a, u$.

Since our proof is based on the local existence result from \cite{danchin2016} for critical initial data, we give the statement of that result here.
\begin{prop} {\rm (}\cite{danchin2016}{\rm )} \label{danchin 2016 crit local thm}
    Let the viscosity coefficients $\lambda$ and $\mu$ depend smoothly on $\rho$ and satisfy \eqref{coeff cndtn}.
    Let $d\geq 2$ and $p \in [1,2d)$, and assume $(a_0,u_0)$ satisfy
    \begin{align*}
        & a_0 \in \db^{ \frac{d}{p} }_{p,1}
        \quad \text{and} \quad 
        u_0 \in \db^{ \frac{d}{p} - 1 }_{p,1},
    \end{align*}
    with, additionally, $1+a_0$ bounded away from $0$.
    Then there exists $T>0$ such that \eqref{CNSo} has a unique local-in-time solution $(a,u)$ with 
    \begin{align*}
        & a\in \widetilde{C}([0,T] ; \db^{ \frac{d}{p} }_{p,1} )
        ,
        \\
        & u \in E_p^{  \frac{d}{p} - 1 }(T)
        \cap 
        \widetilde{L}^\infty_T ( \db^{ \frac{d}{p} - 1 }_{p,1} ).
    \end{align*}
\end{prop}

% The above theorem is proven in Eulerian coordinates only for constant viscosity coefficients with $d\geq 3$ and $p\in [1,d)$, and another proof is given in Lagrangian coordinates for smooth viscosity coefficients and the full range of $d$ and $p$.
% Thanks to the additional regularity assumptions we make, we can prove the full range of $d$ and $p$ using only the Eulerian coordinate proof.

We now give our local existence result.

\begin{prop}{\label{local exist thm hi reg}}
    Let $d\geq 2$, 
    $k\in \{ 1,2,...,d \}$, 
    $p \in [1,2d)$. Then there exists a small constant $c = c(p,d,\mu,\lambda,P)$ and $T > 0$ such that, if
    \begin{align*} 
        \norm{a_0}_{\db^{d/p}_{p,1} \cap \db^{d/p + 1}_{p,1}}
        +
        \norm{x_k a_0}_{\db^{d/p}_{p,1} \cap \db^{d/p + 1}_{p,1}}
        +
        \norm{u_0}_{\db^{d/p - 1}_{p,1} \cap \db^{d/p}_{p,1}}
        +
        \norm{x_k u_0}_{\db^{d/p}_{p,1}} 
        \leq c,
    \end{align*}
    then there exists a unique local solution $(a,u)$ to \eqref{31} and a unique local solution $(A_k,U_k)$ to \eqref{31x} on $[0,T]$, such that

    \begin{align}
        & a \in \widetilde{C}([0,T] ; \db^{ \frac{d}{p} }_{p,1} \cap \db^{\frac{d}{p} + 1}_{p,1} ), \label{loc sp 1}
        \\
        & A_k \in \widetilde{C}([0,T] ; \db^{ \frac{d}{p}  + 1}_{p,1} ), \label{loc sp 2}
        \\
        & u \in E^{ \frac{d}{p} -1 }_p (T) \cap E_p^{ \frac{d}{p} } (T)
        \cap \widetilde{L}^\infty_T ( \db^{ \frac{d}{p} -1  }_{p,1} 
        \cap \db^{ \frac{d}{p}  }_{p,1} 
        ), \label{loc sp 3}
        \\
        & U_k \in E_p^{ \frac{d}{p} } (T)
        \cap \widetilde{L}^\infty_T ( 
        \db^{ \frac{d}{p}  }_{p,1} 
        )
        . \label{loc sp 4}
    \end{align}
\end{prop}

    We prove this result by constructing a series of approximate solutions which are bounded in the function spaces listed above in \eqref{loc sp 1}--\eqref{loc sp 4} and showing that they converge to solutions of \eqref{31} and \eqref{31x} in a suitable function space.
    We first go through this process for $(a,u)$, as the proof for $(A_k,U_k)$ relies on $(a,u)$. 
    The key difference between our construction of $(a,u)$ and Danchin's is the necessity for a higher regularity assumption on $a_0$ and $u_0.$ 
    We will then see in our construction of $(A_k,U_k)$ why this higher regularity was necessary, and why the space $U_k$ lies in $E_p^{ \frac{d}{p} }(T)$ and not $E^{ \frac{d}{p} - 1 }_p(T)$.

\begin{pf}
First, we assume
\begin{align*}
        \norm{a_0}_{\db^{d/p}_{p,1} \cap \db^{d/p + 1}_{p,1}}
        +
        \norm{x_k a_0}_{\db^{d/p}_{p,1} \cap \db^{d/p + 1}_{p,1}}
        +
        \norm{u_0}_{\db^{d/p - 1}_{p,1} \cap \db^{d/p}_{p,1}}
        +
        \norm{x_k u_0}_{\db^{d/p}_{p,1}} 
        \leq c,
    \end{align*}
    for some $c>0.$
Next, we define for all $n \in \mathbb{N} \cup \{0\}$,
\begin{align*}
    & a_0^n \coloneqq \dot{S}_n a_0, 
    \\
    & u_0^n \coloneqq \dot{S}_n u_0,
\end{align*}
which we will use as initial data in our approximate systems. Next, the first term $(a^0,u^0)$ in our sequence of approximations of $(a,u)$ is defined as follows:
\begin{align*}
    & a^0 \coloneqq a^0_0,
    \\
    & u^0 \coloneqq e^{t \mathcal{A} } u^0_0 ,
\end{align*}
where $( e^{t \mathcal{A} } )_{t \geq 0}$ is the semigroup associated with the Lam\'{e} equation \eqref{Lame}.

Then, using Proposition \ref{thm1} to handle $u^0,$ we easily obtain
\begin{align*}
    & a^0 \in 
    C([0,T] ; \db^{ \fdp }_{p,1} \cap \db^{ \fdp + 1}_{p,1} ),
    \\
    & u^0 \in 
    E_p^{ \frac{d}{p} - 1 }(T) \cap E_p^{ \frac{d}{p} }(T),
\end{align*}
for all $T>0.$ 
Next, for any given pair $(a^n,u^n)$ with $n\in \mathbb{N} \cup \{0\}$, we define $(a^{n+1},u^{n+1})$ as the solution of the following system:
\begin{align} \label{32} \displaystyle
    \begin{cases} \displaystyle
        \displaystyle \pt_t a^{n+1} + u^n \cdot \N a^{n+1} = -(1+a^n) \div{u^n} ,
        \\
        \displaystyle \pt_t u^{n+1} - \mathcal{A} u^{n+1} = - u^n \cdot \N u^n - \frac{a^n}{1+a^n} \mathcal{A} u^n - \frac{P'(1+a^n)}{1+a^n} \N a^n,
    \end{cases}
\end{align}
with initial data $(a^{n+1}_0, u^{n+1}_0)$. 
By applying Propositions \ref{thm2} and \ref{thm1}, it can be shown, as in Danchin~\cite{danchin2016}, that for sufficiently small $c$, there exists $T>0$ such that $(a_n)_{n\in \mathbb{N}} $ is a bounded sequence in $C([0,T];\db^{\frac{d}{p}}_{p,1} \cap \db^{\frac{d}{p} + 1}_{p,1} )$ and $(u^n)_{n\in \mathbb{N}}$ is bounded in $E^{ \frac{d}{p} - 1 }_p(T) \cap E_p^{\frac{d}{p}}(T)$. 

Next, we show that these sequences converge strongly in some space by considering, for all $n\in\mathbb{N}$, the differences 
\begin{align*}
    & \delta a^n \coloneqq a^{n+1} - a^n, 
    \\
    & \delta u^n \coloneqq u^{n+1} - u^n, 
\end{align*}
which solve the following system:
\begin{align} \label{41}
    \begin{cases}
    \displaystyle
        \pt_t \da{n+1} + u^{n+1} \cdot \N \da{n+1} = \sum^{3}_{i=1} \delta F^n_i,
        \\
        \displaystyle
        \pt_t \du{n+1} - \mathcal{A} \du{n+1} =
        \sum^{5}_{i=1} \delta G^n_i,
    \end{cases}
\end{align}
where
\begin{align*}
    & \delta F^n_1 \coloneqq - \du{n} \cdot \N a^{n+1}, \quad  \quad 
    \delta F^n_2 \coloneqq - \da{n} \div{u^{n+1}}, \quad  \quad 
    \delta F^n_3 \coloneqq - (1+a^n) \div{\du{n}}, 
    \\
    & \delta G^n_1 \coloneqq 
        \bigc{ \frac{a^n}{1+a^n} - \frac{a^{n+1}}{1+a^{n+1}} } \mathcal{A} u^{n+1}, \quad  \quad 
    \delta G^n_2 \coloneqq 
        - \frac{a^n}{1+a^n} \mathcal{A} \du{n},
    \\
    & \delta G^n_3 \coloneqq
        \frac{P'(1+a^n)}{1+a^n} \N a^n 
         - \frac{P'(1+a^{n+1})}{1+a^{n+1}} \N a^{n+1}, 
    \\
    & 
    \delta G^n_4 \coloneqq
        -u^{n+1} \cdot \N \du{n}, 
        \quad  \quad 
    \delta G^n_5 \coloneqq
        - \du{n} \cdot \N u^n .
\end{align*}
We consider $(\delta a^n, \delta u^n)$ in the space 
$C([0,T] ; \db^{d/p}_{p,1}  ) \times E_p^{ \frac{d}{p} - 1 }(T)$. 
First, we check the regularity of $(\da{0},\du{0})$.
Since $\da{0} \coloneqq a^1 - a^0,$ and $a^0 \coloneqq \dot{S}_0 a_0$, we can add the term $u^0 \cdot \N a^0$ to both sides of the equation for $a^1$ to obtain the following equation for $\da{0}$:
\begin{align*}
    \pt_t \da{0} + u^0 \cdot \N \da{0} = -(1+a^0) \div{u^0} -u^0 \cdot \N a^0.
\end{align*}
Then, by Proposition \ref{thm2}, 
\begin{align*}
    & \norm{\da{0}}_{ \widetilde{L}^\infty_T \db^{d/p}_{p,1} }
    \leq
    e^{C \int^T_0 \norm{ \N u^0 }_{\db^{d/p}_{p,1}} \dd t}
    \biggc{
    \norm{\da{0}}_{\db^{d/p}_{p,1}} + \int^T_0 \norm{ 
    (1+a^0) \div{u^0} + u^0 \cdot \N a^0
    }_{\db^{d/p}_{p,1}} \dd t
    }
    \\
    & \leq C
    \biggc{
    \norm{ a_0 }_{\db^{d/p}_{p,1}}
    +
    \norm{u^0}_{L^1_T \db^{d/p +1}_{p,1}}
    + 
    \norm{a^0}_{L^\infty_{T} \db^{d/p}_{p,1} }
    \norm{u^0}_{L^1_T \db^{d/p +1}_{p,1}}
    +
    \norm{a^0}_{L^\infty_{T} \db^{d/p + 1}_{p,1} }
    \norm{u^0}_{L^1_T \db^{d/p}_{p,1}}
    }
    \\
    & \leq
    C(1+T)c.
\end{align*}
Next, since $u^0 \coloneqq e^{t\mathcal{A}} \dot{S}_0 u_0$, $u^0$ solves 
\begin{align*}
    \begin{cases}
        \pt_t u^0 - \mathcal{A} u^0 = 0, \text{ for } t>0, \ x\in\R^d,
        \\
        u^0(0) = \dot{S}_0 u_0, \text{ for } x\in\R^d.
    \end{cases}
\end{align*}
Combining with the equation for $u^1,$ we get
\begin{align*}
    \begin{cases}
    \displaystyle
        \pt_t \du{0} - \mathcal{A} \du{0} =
            -u^0 \cdot\N u^0 - \frac{a^0}{1+a^0} \mathcal{A} u^0
            -\frac{P'(1+a^0)}{1+a^0} \N a^0,
            \text{ for } t>0, \ x \in \R^d
        \\
        \displaystyle
        \du{0}(0) = \dot{\Delta}_1 u_0, \text{ for } x\in\R^d,
    \end{cases}
\end{align*}
The right-hand side of the above equation is the same as the equation for $u^1$, so we know that applying Proposition \ref{thm1} gives
\begin{align*}
    \norm{\du{0}}_{E_p^{ \frac{d}{p} - 1 }(T)} < C 
    \biggc{
    \norm{ u_0 }_{\db^{ \frac{d}{p} - 1 }_{p,1}}
    +
    cT
    },
\end{align*}
for sufficiently small $c$. So we have $(\da{0},\du{0}) \in L^\infty_T \db^{ \frac{d}{p} }_{p,1} \times E_p^{ \frac{d}{p} }(T) $. Next, let $n\in \mathbb{N} \cup \{ 0 \}$. Using Propositions \ref{thm2} and \ref{thm1}, and the boundedness of 
$(a_n)_{n\in \mathbb{N}} $ in $C([0,T];\db^{\fdp}_{p,1} \cap \db^{ \fdp + 1}_{p,1} )$ 
and 
$(u^n)_{n\in \mathbb{N}}$ in $E_p(T) \cap E_p^{+1}(T)$, we obtain the following bounds, again taking $c$ sufficiently small:
\begin{align*}
    & \norm{\da{n+1}}_{L^\infty_T \db^{\fdp}_{p,1}}
        \leq C \norm{ \da{n+1} (0) }_{\db^{\fdp}_{p,1}}
        + \frac{1}{8} \norm{ \da{n} }_{ L^\infty_T \db^{ \fdp }_{p,1} }
        + 2 \norm{ \du{n} }_{E_p(T)}
    \\
    & \norm{ \du{n+1} }_{E_p(T)}
        \leq C \norm{ \du{n+1}(0) }_{\db^{ \fdp - 1}_{p,1}}
        + \frac{1}{8} \biggc{
            \norm{\da{n}}_{L^\infty_T \db^{ \fdp }_{p,1}}
            + \norm{ \du{n} }_{E_p(T)}.
        }
\end{align*}
Combining both, we have 
\begin{align*}
    & \norm{\da{n+1}}_{L^\infty_T \db^{d/p}_{p,1}}
    + 4 \norm{ \du{n+1} }_{E_p(T)}
    \\
    & \quad \leq 
    C \biggc{
    \norm{ \da{n+1} (0) }_{\db^{d/p}_{p,1}}
    + \norm{ \du{n+1}(0) }_{\db^{d/p - 1}_{p,1}}
    }
    + \frac{5}{8}
    \biggc{
    \norm{ \da{n} }_{ L^\infty_T \db^{d/p}_{p,1}}
    + 4 \norm{ \du{n} }_{E_p(T)}
    }.
\end{align*}
We recall that, for a pair of sequences $(l_n)_{n\in\mathbb{N}}$, $(L_n)_{n\in\mathbb{N}}$, if 
\begin{align*}
    & \sum^{\infty}_{n=0} l_n < \infty,
    \\
    & K < 1,
    \\
    & L_{n+1} \leq C l_n + K L_n,
\end{align*}
then
\begin{align*}
    \sum^{\infty}_{n=0} L_n < \infty.
\end{align*}
We thus conclude that
\begin{align*}
    & 
    \sum^{\infty}_{n=0} \norm{\da{n+1}}_{L^\infty_T \db^{d/p}_{p,1}}
    + 4 \norm{ \du{n+1} }_{E_p(T)} <\infty,
\end{align*}
and so, there exist functions $a \in L^\infty_T (\db^{ \frac{d}{p} }_{p,1}) $ and $u \in E_p^{ \frac{d}{p} - 1 }(T) $ such that
\begin{align*}
    & a^n \to a \text{ in } L^\infty_T (\db^{ \frac{d}{p} }_{p,1}),
    \\
    & u^n \to u \text{ in } E_p^{ \frac{d}{p} - 1 }(T).
\end{align*}
Next, it can be checked that each term in \eqref{32} converges to their respective corresponding term in \eqref{31} in the strong topology of $L^1_T (\db^{ \frac{d}{p}  - 1}_{p,1})$. Thus we conclude that $(a,u)$ is a strong solution of \eqref{31}.
% {\bl (Actually, should I only use the abbreviation $L^1_T (\db^{ \frac{d}{p}  - 1}_{p,1})$ when writing norms? Should I stick with the full notation $L^1 ( 0,T; \db^{ \frac{d}{p}  - 1}_{p,1})$ when writing the set?)}

Also, by the boundedness of $(a_n)_{n\in\mathbb{N}}$ and $(u_n)_{n\in\mathbb{N}}$, and by the Banach-Alaoglu Theorem, we get weak-$\ast$ convergence to the same functions in the subcritical spaces:
\begin{align*}
    & a^n \weakstar a \text{ in } L^\infty_T (\db^{ \frac{d}{p} + 1 }_{p,1}),
    \\
    & u^n \weakstar u \text{ in } L^\infty_T (\db^{ \frac{d}{p} }_{p,1}).
\end{align*}

Uniqueness of the solution may be checked exactly as in \cite{danchin2016}. In particular, for two solutions $(a^1,u^1)$ and $(a^2,u^2)$ of \eqref{31}, it can be shown that
\begin{align*}
    (a^2 - a^1 , u^2 - u^1) \equiv 0 \text{ in } C([0,T] ; \db^{\fdp - 1}_{p,1}) \times F_p(T),
\end{align*}
where
\begin{align*}
    F_p(T) \coloneqq C([0,T] ; \db^{\fdp - 2}_{p,1}) 
    \cap L^1_T (\db^{\fdp}_{p,1}).
\end{align*}
% {\bl (This is the `one derivative lower' space that Danchin uses for convergence. It should be possible to directly prove uniqueness in a higher regularity space, using my extra regularity, but that isn't necessary.)}

We now move onto $(A_k,U_k)$. Since these functions are defined in terms of $(a,u)$, and we already have the above results proven for $(a,u)$, it will suffice only to construct bounded (and thus weakly-$\ast$ convergent) sequences $(A_{k,n})_{n\in \mathbb{N}}$ and $(U_{k,n})_{n \in \mathbb{N}}$ in the correct spaces.

We follow steps similar to before, first defining
\begin{align*}
    & A^n_{k,0} \coloneqq x_k a^n_0 = x_k (\dot{S}_n a_0), 
    \\
    & U^n_{k,0} \coloneqq x_k u^n_0 = x_k (\dot{S}_n u_0),
\end{align*}
and
\begin{align*}
    & A^0_k \coloneqq x_k a^0 = x_k a^0_0,
    \\
    & U^0_k \coloneqq x_k u^0 = x_k (e^{t\mathcal{A}} u^0_0).
\end{align*}
We clearly have 
\begin{align*}
    A^0_k \in C([0,T] ;  \db^{\fdp  + 1}_{p,1} ),
\end{align*}
but $U^0_k$ is a bit less straightforward. By definition of $u^0$, we know that
\begin{align*}
    x_k \biggc{
    \pt_t u^0 - \mathcal{A} u^0 =0
    }.
\end{align*}
Let $e_k$ denote the unit vector along the $x_k$-axis. Rewriting in terms of $U^0_k$ gives
\begin{align*}
    \pt_t U^0_k -\mathcal{A} U^0_k 
    & = -2\mu \pt_k u^0 - (\lambda + \mu) \biggc{
    \div(u^0) e_k - \N (u^0_k)
    }.
    % \\
    % & \eqqcolon
    % S^0.
\end{align*}
Then, by Proposition \ref{thm1}, we have for all $\sigma\in\R$:
\begin{align*}
    & \norm{U^0_k}_{L^\infty_T \db^{\sigma}_{p,1}} 
    + \int^T_0 \norm{ \N^2 U^0_k }_{\db^\sigma_{p,1}} \dd t
    \\
    &
    \leq
    C \biggc{
    \norm{ U^0_k (0) }_{\db^\sigma_{p,1}}
    + \int^T_0 \bignorm{ 
    %S^0(t)
    2\mu \pt_k u^0(t) + (\lambda + \mu) \biggc{
    \div(u^0 (t) ) e_k - \N (u^0_k (t) )
    }
    }_{\db^\sigma_{p,1}} \dd t
    }
    \\
    & \leq 
    C \biggc{
    \norm{ u_0 }_{\db^\sigma_{p,1}}
    + \norm{ U_{k,0} }_{\db^\sigma_{p,1}}
    + \norm{ u^0 }_{L^1_T \db^{\sigma + 1}_{p,1}}
    }.
\end{align*}
We will see in the next step why it is desirable to choose $\sigma = d/p$ instead of $d/p - 1$.

We have thus far defined $A^n_{k,0}, \ U^n_{k,0}, \ A^0_, \text{ and } U^0_k$ simply by multiplying by $x_k$ the corresponding functions related to $a$ and $u$. We continue this pattern now by multiplying the whole approximate system \eqref{32} by $x_k$ to obtain
\begin{align} \label{32x} \displaystyle
    \begin{cases} \displaystyle
        \pt_t A^{n+1}_k + u^n \cdot \N A^{n+1}_k
            = -(x_k + A^n_k) \div{u^n}
            + u^n_k a^{n+1},
        \\ \displaystyle
        \pt_t U^{n+1}_k - \mathcal{A} U^{n+1}_k
            = -U^n_k \cdot \N u^n 
            -\frac{a^n}{1+a^n}\mathcal{A}U^n_k
            - \frac{P'(1+a^n)}{1+a^n} \biggc{
                \N A^n_k - a^n e_k
            }
            \\
            \quad \quad \quad \quad \quad \quad \quad \quad \quad 
            \displaystyle 
            + \frac{ a^n }{1+a^n} B_n - B_{n+1},
    \end{cases}
\end{align}
where 
\begin{align*}
    B_n \coloneqq 
            2\mu\pt_k u^n + (\lambda + \mu) \biggc{
                \div(u^n) e_k + \N (u^n_k)
            }
            .
\end{align*}
Then, by Propositions \ref{thm2} and \ref{thm1}, taking $c$ smaller if necessary, there exists a large constant $C>0$ such that if $T<1$ and 
\begin{align*}
    & \norm{A^n_k}_{L^\infty_T  \db^{\fdp +1}_{p,1} } \leq C(1+T),
    \\
    &
    \norm{U^n_k}_{E^{ \frac{d}{p} }_p(T)} \leq C(1+T),
\end{align*}
then the same inequalities are satisfied by $A^{n+1}_k$ and $U^{n+1}_k$. 
The reason for taking a higher regularity exponent for $U^n_k$ than $u^n$ is the term 
\begin{align} \label{reason hi reg Uk}
\frac{P'(1+a^n)}{1+a^n} a^n e_k.
\end{align}
Since we don't have $a_0 \in \db^{ \fdp - 1}_{p,1}$, we are forced to look at $U^n_k$ in $E^{ \fdp }_p(T)$. 
But then we also need 
\begin{align*}
    \frac{P'(1+a^n)}{1+a^n} \N A^n_k \in L^\infty_T \db^{\fdp }_{p,1},
\end{align*}
hence the requirement that $(A^n_k)_n$ be bounded in $L^\infty_T (\db^{\fdp + 1}_{p,1})$, `one derivative higher than' the critical space for $(a^n)_n$.
Subsequently, we need the additional higher regularity exponents for $a$ and $u$, but this fortunately does not trigger any chain reaction past this point.

We now have that $(A^n_k)_{n\in\mathbb{N}}$ is a bounded sequence in $L^\infty_T ( \db^{\fdp + 1}_{p,1})$ and $(U^n_k)_{n\in\mathbb{N}}$ is a bounded sequence in $E^{\fdp}_p(T)$. 
Thus, by the Banach-Alaoglu Theorem, there exist functions $A_k$ and $U_k$ such that
\begin{align*}
    & A^n_k \weakstar A_k \text{ in } L^\infty_T ( \db^{\fdp + 1}_{p,1}),
    \\
    & U^n_k \weakstar U_k \text{ in } L^\infty_T (\db^{\fdp }_{p,1}).
\end{align*}
It can also be checked that $U_k \in E^{\fdp}_p(T)$.
Combining with our convergence results on $(a^n)_{n\in\mathbb{N}}$
and $(u^n)_{n\in\mathbb{N}}$ confirms that these functions satisfy $A_k = x_k a$ and $U_k = x_k u$, which completes the proof.
\end{pf}

\section{Global Existence}
\label{section global}

For the global existence problem, we rearrange our system \eqref{31} as follows
\begin{align} \label{78}
    \begin{cases} 
        \displaystyle 
        \pt_t a + \div(u) = f \coloneqq -\div(au)
        \\
        \displaystyle 
        \pt_t u - \mathcal{A} u + P'(1)\N a = g \coloneqq 
        -u\cdot\N u - \frac{a}{1+a} \mathcal{A} u - \beta(a) \N a,
    \end{cases}
\end{align}
where
\[
\beta(a) \coloneqq  \frac{P'(1+a)}{1+a} - P'(1).
\]
We next split the above system into its divergence-free and curl-free parts using the projection operators $\P$ and $\Q$. Relabelling $P'(1) \eqqcolon \al$ and $(\lambda + 2\mu) \eqqcolon \nu$, this gives us
\begin{align} \label{79}
    \begin{cases}
        \pt_t a + \div(\Q u) = f,
        \\
        \pt_t \Q u - \nu \Delta \Q u + \al \N a = \Q g,
        \\
        \pt_t \P u - \mu \Delta \P u = \P g.
    \end{cases}
\end{align}
Looking at the bottom line, we see that $\P u$, the incompressible part of the velocity, obeys a heat equation.
The top two lines, governing the density and curl-free part of the velocity, may be analysed separately as a $2\times 2$ system in terms of the following scalar function
\begin{align*}
    v \coloneqq |D|^{-1} \div(u),
\end{align*}
Since the difference between $v$ and $\Q u$ is only a $0$-order Fourier multiplier, bounding $v$ in any homogeneous Besov norm is equivalent to bounding $\Q u$, by Proposition \ref{fm est for besov}.
See \cites{danchin2016, daithi2023} for further discussion.

Next, we introduce the following rescaling:
\[
a(t,x) = \tilde{a} \biggc{ \frac{\al}{\nu} t, \frac{\sqrt{\al}}{\nu} x },
\quad
u(t,x) = \sqrt{\al} \tilde{u} \biggc{ \frac{\al}{\nu} t, \frac{\sqrt{\al}}{\nu} x }
\]
and observe that $(\tilde{a} , \tilde{u} )$ solve \eqref{79} with $\al = \nu = 1$. 
Thus, we may assume without loss of generality that $\al = \nu = 1$.

% %%%%%%%%old lines about v
% We then have that $(a,v)$ solves the following system:
% \begin{align}
%     \begin{cases}
%         \pt_t a + |D|v = f,
%         \\
%         \pt_t v - \Delta v - |D| a = h \coloneqq |D|^{-1} \div(g).
%     \end{cases}
% \end{align}
% We then note that the Fourier transform of $(a,v)$ solves the following vector equation:
% \begin{align}
%     \pt_t 
%     \begin{bmatrix}
%         \ha \\ \hv    
%     \end{bmatrix}
%     =
%     M_{|\x|} 
%     \begin{bmatrix}
%         \ha \\ \hv    
%     \end{bmatrix}
%     +
%     \begin{bmatrix}
%         \hat{f} \\ \hat{h}    
%     \end{bmatrix},
%     \quad 
%     \text{where}
%     \quad
%     M_{|\x|} 
%     \coloneqq
%     \begin{bmatrix}
%         0 \ -|\x|
%         \\
%         |\x| \ -|\x|^2
%     \end{bmatrix}.
% \end{align}
% By analysing the homogeneous case (as in \cite{daithi2023}), where we set $f=h=0$, it can be easily shown that the linear solution 
% \[
% e^{tM} \begin{bmatrix}
%         a(0) \\ v(0)    
%     \end{bmatrix}
%     \coloneqq
%     \mathcal{F}^{-1}
%     \Bigg{[}
% e^{tM_{|\x|}} \begin{bmatrix}
%         \ha(0) \\ \hv(0)    
%     \end{bmatrix}
%     \Bigg{]}
% \]
% decays in $L^2$ like the fundamental solution to a heat equation.
% That is,
% \[
% \bignorm{
% e^{ tM_{|\x|} } \begin{bmatrix}
%         \ha(0) \\ \hv(0)    
%     \end{bmatrix}
% }_2
% \leq C(t+1)^{- d/4 }.
% \]
% In fact, the above linear solution decays even faster than heat for $p>2$, but we will only need boundedness in $L^2$ for our purposes here.

We can give the exact same treatment to $(A_k,U_k)$. Starting of by rewriting \eqref{31x} as
\begin{align} \label{78x}
    \begin{cases}
        \pt_t A_k + \div(U_k) = \mathfrak{f} 
        \coloneqq 
        -A_k \div(u) - U_k \cdot \N a + u_k,
        \\
        \pt_t U_k - \mathcal{A} U_k + P'(1) \N A_k = \mathfrak{g}
        \coloneqq 
        -2 \mu \pt_k u + a e_k 
        \\
         \quad \quad \quad \quad 
        -(\lambda + \mu) \biggc{
        \div(u) e_k +\N (u_k)
        }
        - (x_k \beta) \N a
        \displaystyle -U_k \cdot \N u -\frac{A_k}{1+a} \mathcal{A}u.
    \end{cases}
\end{align}
Splitting into divergence-free and curl-free parts, and using the same notation as before, this becomes
\begin{align} \label{79x}
    \begin{cases}
        \pt_t A_k + \div(\Q U_k) = \mathfrak{f},
        \\
        \pt_t \Q U_k - \nu \Delta \Q U_k + \al \N A_k = \Q \mathfrak{g},
        \\
        \pt_t \P U_k - \mu \Delta \P U_k = \P \mathfrak{g}.
    \end{cases}
\end{align}

% %old argument about $V_k$
% Then making the same rescaling argument and setting
% \[
% V_k \coloneq |D|^{-1} \div(U_k)
% \]
% leads to this simplified system for the curl-free part:
% \begin{align}
%     \begin{cases}
%         \pt_t A_k + |D|V_k = \mathfrak{f},
%         \\
%         \pt_t V_k - \Delta V_k - |D| A_k = \mathfrak{h} \coloneqq |D|^{-1} \div(\mathfrak{g}).
%     \end{cases}
% \end{align}

Crucially, the left hand side of the above systems are the same for $(a,u)$ and $(A_k,U_k)$, and so the same analysis holds for their linear solutions. 
It is the nonlinear part of solutions where they differ.
Most important is the presence of terms such as $u_k$ in $\mathfrak{f}$ and $a e_k$ in $\mathfrak{g}$.
These terms prevent us from bounding the low-frequency norm of $(A_k,U_k)$ with the same regularity exponent as $(a,u)$.
Indeed, for the low-frequency part of $(a,u)$ and $(A_k,U_k)$, we apply the following lemma.
\begin{lem}{\rm (}\cite{danchin2016}{\rm )} \label{danchin prop 7}
    Let $s\in\R$ and $(a,u)$ solve \eqref{78} (or equivalently $(A_k,U_k)$ solve \eqref{78x}) with $P'(1) = \nu = 1$. Let $j_0\in\Z$ be the frequency cut-off constant.
    Then there exists a constant $C=C(j_0,\mu)$ such that for $t > 0$,
    \begin{align}
        & \norm{ (a,u) }^l_{\widetilde{L}^\infty_t \db^s_{2,1} }
        +
        \norm{ (a,u) }^l_{ L^1_t \db^{s+2}_{2,1} }
        \leq
        C
        \biggc{
        \norm{ (a_0,u_0) }^l_{\db^s_{2,1} }
        +
        \norm{ (f,g) }^l_{ L^1_t \db^{s}_{2,1} }
        }, \label{lo-freq est prop7}
        \\
        & \norm{ a }^h_{\widetilde{L}^\infty_t \db^{s+1}_{2,1} }
        +
        \norm{ a }^h_{ L^1_t \db^{s+1}_{2,1} }
        +
        \norm{ u }^h_{\widetilde{L}^\infty_t \db^{s}_{2,1} }
        +
        \norm{ u }^h_{ L^1_t \db^{s+2}_{2,1} }
        \notag
        \\
        & \quad \leq
        C
        \biggc{
        \norm{ a_0 }^h_{\db^{s+1}_{2,1} }
        +
        \norm{ u_0 }^h_{\db^{s}_{2,1} }
        +
        \norm{ f }^h_{L^1_t \db^{s+1}_{2,1} }
        +
        \norm{ g }^h_{L^1_t \db^{s}_{2,1} }
        },
    \end{align}
    and similar for $(A_k,U_k).$
\end{lem}

For our global existence result, we will use inequality \eqref{lo-freq est prop7} with $s=d/2 - 1$ to bound the low frequencies of $(a,u)$. 
We see that, while on the left-hand side we have a regularity exponent of $d/2 + 1$ in the $ L^1_t$-norm, we only have a regularity exponent of $d/2 -1$ on the right-hand side $ L^1_t$-norm of $(f,g)$.
This gap in the exponents is handled by the presence of derivatives in every term of the original nonlinear terms $f$ and $g$ and by product estimates.
However, since, for example, the weighted nonlinear term $\mathfrak{f}$ has the term $u_k$, we cannot hope to use \eqref{lo-freq est prop7} for $(A_k,U_k)$ even with the higher regularity exponent $s=d/2$, unless we can improve the low-frequency regularity of $(a,u)$ first. 
This will be achieved by applying a low-frequency decay estimate on $(a,u)$ due to Danchin-Xu~(\cite{Danchin-Xu2017}).

The high-frequency norms, on the other hand, can be bounded by the same methods as in \cite{danchin2016}. 
The inequalities found therein are readily extended to the subcritical norms for $(a,u)$ and also readily applied to $(A_k,U_k)$. 

We first give a global existence result for $(a,u)$ in the combined critical-and-subcritical regime. 
Only afterwards will we prove our final result which includes the weighted Besov regime.
For this result, we introduce the function space $Y_p$ as the set of all pairs of functions $(a,u)$, where $a : [0,\infty) \times \R^d \to [0,\infty)$ is a scalar function and $u: [0,\infty) \times \R^d \to \R^d$ is a $d$-vector function, satisfying the following:
\begin{align*}
    & (a,u)^l \in \widetilde{C} (\R_{>0}; \db^{\frac{d}{2} - 1}_{2,1}) \cap L^1(\R_{>0}; \db^{\frac{d}{2} + 1}_{2,1}), 
    \\
    & a^h \in \widetilde{C} (\R_{>0}; \db^{\frac{d}{p} + 1}_{p,1}) \cap L^1(\R_{>0}; \db^{\frac{d}{p} + 1}_{p,1}),
    \\
    & u^h \in \widetilde{C}(\R_{>0}; \db^{\frac{d}{p}}_{p,1})
    \cap L^1(\R_{>0}; \db^{\frac{d}{p} + 2}_{p,1})
    % ,
    % \\
    % & (x_k a)^l \in \widetilde{C} (\R_{>0}; \db^{\frac{d}{2}}_{2,1}) 
    % \cap L^1_t ( \R_{>0} ; \db^{\frac{d}{2} + 2}_{2,1} ), 
    % \quad 
    % (x_k a)^h \in \widetilde{C} (\R_{>0}; \db^{\frac{d}{p} +1 }_{p,1})
    % \cap L^1_t ( \R_{>0} ; \db^{\frac{d}{p} + 1}_{p,1} ),
    % \\
    % & (x_k u)^l \in \widetilde{C} (\R_{>0}; \db^{\frac{d}{2}}_{2,1})
    % \cap L^1_t ( \R_{>0} ; \db^{\frac{d}{2} + 2}_{2,1} ),
    % \quad
    % (x_k u)^h \in \widetilde{C} (\R_{>0}; \db^{\frac{d}{p}}_{p,1})
    % \cap L^1_t ( \R_{>0} ; \db^{\frac{d}{p} + 2}_{p,1} )
    .
\end{align*}
$Y_p$ is then equipped with the obvious norm corresponding to the strong topologies for the above spaces.
In words, the definition of $Y_p$ is the same as $S_p$ from our main theorem, but without the regularity in weighted Besov spaces. 
%{\bl consider rewriting}

\begin{prop} \label{gl ext Danchin extended} 
%{\bl (Global existence extended from Danchin to include hi derivs)}
    Let $d\geq2$ and $p \in [2, \text{min} \{4, 2d/(d-2)\} ]$, with $p\neq 4$ in the $d=2$ case.
    Assume 
    without loss of generality 
    that $P'(1) = \nu = 1$.
    Then there exist a frequency cut-off constant $j_0 \in \Z$ and a small constant $c = c(p,d,\mu,P) >0$ such that, if $(a_0,u_0)$ 
    %and $(A_{k0},U_{k0})$ 
    satisfy 
    \begin{align}
    \notag
        Y_{p,0} \coloneqq &  
        \norm{(a_0,u_0)}^l_{\dot{B}^{\frac{d}{2} - 1}_{2,1}} 
        + \norm{a_0}^h_{\db^{\frac{d}{p} + 1}_{p,1}} 
        + \norm{u_0}^h_{\db^{\frac{d}{p}}_{p,1}}
        \leq c,
    \end{align}
    then \eqref{31} has a unique solution $(a,u)$ in the space $Y_p$. 
    Also, there exists a constant $C = C(p,d,\mu, P, j_0)$ such that
    \begin{align} \label{subcrit est}
        \norm{(a,u)}_{Y_p} \leq C Y_{p,0}.
    \end{align}
\end{prop}

% The proof is identical that of the main theorem in \cite{danchin2016}, but with a Besov index that is $1$ higher for the high-frequency norms. 
% This is of course possible thanks to our higher regularity assumption on the high frequencies of the initial data.

The proof rests on an a priori bound for $(a,u)$.
We give an outline of how this a priori bound is obtained, which uses the same steps as in \cite{danchin2016}, but with a higher regularity exponent for the high-frequency norm.
These same steps can then be adapted for the global existence proof in weighted Besov spaces.
% For a more detailed proof with more exposition on the individual steps, see \cite{danchin2016}.
We define the following quantity:
\begin{align*}
    Y_p(t) & \coloneqq 
    \norm{ (a,u) }^l_{\widetilde{L}^\infty_t \db^{ \frac{d}{2} - 1}_{2,1} }
    + 
    \norm{ (a,u) }^l_{ L^1_t \db^{ \frac{d}{2} + 1}_{2,1} }
    +
    \norm{a}^h_{ \widetilde{L}^\infty_t \db^{ \fdp + 1 }_{p,1} }
    +
    \norm{a}^h_{L^1_t \db^{ \fdp + 1 }_{p,1} }
    \notag
    \\
    & \quad +
    \norm{ u }^h_{ \widetilde{L}^\infty_t \db^{ \fdp }_{p,1} }
    +
    \norm{ u }^h_{L^1_t \db^{ \fdp + 2 }_{p,1} }.
\end{align*}
Our goal is to show that $Y_p(t)$ satisfies \eqref{subcrit est} for all $t>0$.
Throughout the proof of the a priori bound, we consider a smooth solution $(a,u)$ such that $\norm{a}_{L^\infty L^\infty} \leq 1/2$.

The low frequencies are first bounded by Lemma \ref{danchin prop 7} with $s= d/2 - 1$, giving
\begin{align} \label{gl lin est for low freqs}
    \norm{ (a,u) }^l_{\widetilde{L}^\infty_t \db^{\frac{d}{2} - 1}_{2,1}}
    +
    \norm{ (a,u) }^l_{ L^1_t \db^{\frac{d}{2}  + 1}_{2,1}}
    \leq C
    \biggc{
    \norm{ ( a_0,u_0 ) }^l_{ \db^{\frac{d}{2}  - 1}_{2,1}}
    +
    \norm{ ( f,g ) }^l_{ L^1_t \db^{ \frac{d}{2}  -1 }_{2,1}}
    }.
\end{align}

For the high frequencies, recalling \eqref{78}, we use the fact that $\P u$ solves a heat equation and apply Proposition \ref{thm1} with $s=d/p$ (one higher than the proof in the critical regime) to get
\begin{align} \label{gl lin est for Pu}
    \norm{ \P u }^h_{\widetilde{L}^\infty_t \db^{\frac{d}{p}}_{p,1} }
    +
    \mu \norm{ \P u }^h_{ L^1_t \db^{\frac{d}{p} + 2}_{p,1} }
    \leq
    C \biggc{
    \norm{ \P u_0 }^h_{  \db^{\frac{d}{p}}_{p,1} }
    +
    \norm{ \P g }^h_{ L^1_t \db^{\frac{d}{p} }_{p,1} }
    }.
\end{align}

Next, we define the `effective velocity' 
\begin{align*}
    w \coloneqq \N (-\Delta)^{-1} (a - \div(u))
\end{align*}
and find that it also solves a heat equation:
\begin{align*}
    \pt_t w - \Delta w =
    \N (-\Delta)^{-1} (f - \div(g)) + w - (-\Delta)^{-1}\N a.
\end{align*}
Thus we may also apply Thereom \ref{thm1} to $w$ and get
\begin{align} \label{gl lin est for w}
    & \notag 
    \norm{ w }^h_{ \widetilde{E}^{ \fdp }_p (t) } 
    \coloneqq
    \norm{ w }^h_{\widetilde{L}^\infty_t \db^{\fdp }_{p,1} }
    +
    \norm{ w }^h_{ L^1_t \db^{\fdp  + 2}_{p,1} }
    \\
    & \leq C \biggc{
    \norm{w_0}^h_{\db^{ \frac{d}{p} }_{p,1}}
    +
    \norm{ f - \div(g) }^h_{ L^1_t \db^{\frac{d}{p} - 1 }_{p,1}} 
    +
    \norm{ w - (-\Delta)^{-1} \N a }^h_{ L^1_t \db^{\frac{d}{p} }_{p,1}} 
    }
\end{align}
Finally, noting that
\[
\pt_t a + \div(au) + a = -\div(w),
\]
we may obtain by steps similar to the proof of Proposition \ref{thm2},
\begin{align} \label{gl lin est for a}
    & \notag
    \norm{ a }^{h}_{\widetilde{L}^\infty_t \db^{\fdp +1}_{p,1}}
    +
    \norm{ a }^{h}_{ L^1_t \db^{\fdp +1}_{p,1} }
    %\\
    %& 
    \leq C \notag
    \Bigc{
    \norm{ a_0 }^h_{ \db^{ \fdp  + 1 }_{p,1} }
    +
    \norm{ \div(w) }^h_{ L^1_t \db^{\fdp +1}_{p,1} }
    \\
    & \quad \quad \quad \quad +
    \int^{t}_{0} (
    \norm{\N u}_{  \db^{ \fdp  }_{p,1} }
    \norm{a}_{ \db^{ \fdp  + 1 }_{p,1} }
    +
    \norm{\N u}_{  \db^{ \fdp  + 1 }_{p,1} }
    \norm{a}_{  \db^{ \fdp  }_{p,1} }
    )
    \dd \tau
    }.
\end{align}
Observing that
\begin{align*}
    & \norm{ u }^h_{ \widetilde{E}^{\fdp }_p (t) } 
    % \coloneqq
    % \norm{ u }^h_{\widetilde{L}^\infty_t \db^{d/p}_{p,1} }
    % +
    % \norm{ u }^h_{ L^1_t \db^{d/p + 2}_{p,1} }
    %\\
    %& 
    \leq 
    \norm{ \P u }^h_{ \widetilde{E}^{\fdp }_p (t) }
    +
    \norm{ w }^h_{ \widetilde{E}^{\fdp }_p (t) }
    +
    C 
    \Bigc{
    \norm{a}^h_{ \widetilde{L}^\infty_t \db^{\fdp  - 1}_{p,1}}
    +
    \norm{a}^h_{ L^1_t \db^{\fdp  + 1}_{p,1}}
    },
\end{align*}
we may combine \eqref{gl lin est for low freqs}, \eqref{gl lin est for Pu}, \eqref{gl lin est for w}, and \eqref{gl lin est for a} (for high enough frequency cut-off constant $j_0$) to arrive at the following
\begin{align*}
    Y_p(t) \leq \notag
    &
    C 
    \Biggc{
    Y_{p,0}
    +
    \int^t_0
    (
    \norm{ (f,g) }^l_{ \db^{\frac{d}{2} - 1}_{2,1} }
    +
    \norm{ f }^h_{ \db^{ \fdp - 1 }_{p,1} }
    +
    \norm{ g }^h_{ \db^{ \fdp }_{p,1} }
    \\ &
    \quad \quad \quad \quad 
    +
    \norm{\N u}_{  \db^{ \fdp   }_{p,1} }
    \norm{a}_{ \db^{ \fdp  + 1 }_{p,1} }
    +
    \norm{\N u}_{  \db^{ \fdp  + 1 }_{p,1} }
    \norm{a}_{  \db^{ \fdp }_{p,1} }
    )
    \dd \tau
    }.
\end{align*}

The final step is to show that the integral above is less than $C Y_p(t)^2$. 
This is shown using simple product laws and, for the low frequencies of $f$ and $g$, estimates on their bony decompositions.
See \cite{danchin2016} and \cite{danchin-he2016} for details.

We then obtain that, for all $t>0$,
\begin{align*}
    Y_p(t) \leq C \bigc{
    Y_{p,0} + Y_p(t)^2
    }.
\end{align*}
Then, as long as
\begin{align*}
    2 C Y_p (t) \leq 1,
\end{align*}
we have
\begin{align} \label{supercrit a priori est}
    Y_p(t) \leq 2CY_{p,0}.
\end{align}

\noindent\textbf{Proof of Proposition \ref{gl ext Danchin extended}.} 
Under the assumptions of Proposition \ref{gl ext Danchin extended}, if $Y_{p,0}$ is sufficiently small, then by Proposition \ref{local exist thm hi reg} (without the conditions on weighted Besov norms), we have unique local-in-time existence of $(a,u)$ on the time interval $[0,T]$, for some $T>0$.
Next, suppose that $T^* < \infty$ is the maximal existence time 
for $(a,u)$.
However, again if $Y_{p,0}$ is taken small enough, our steps above obtaining the a priori estimate \eqref{supercrit a priori est} show that, for all $t<T^*$,
\[
Y_p(t) \leq C Y_{p,0}.
\]
Then if $C Y_{p,0} \leq c$, where $c\leq 1/2$ is the constant from Proposition \ref{local exist thm hi reg}, we may choose $t_0 \in [0,T^*]$, treat $(a(t_0),u(t_0))$ as initial data, and reapply Proposition \ref{local exist thm hi reg} to show existence on a new time interval $[t_0, t_0+T]$, where $T$ is independent of $t_0$.
We may choose $t_0 > T^* - T$ to extend our solution past $T^*$, but this contradicts our maximality assumption on $T^*$.
\hfill\qed

\bigskip
%Now we move on to our last proposition, where we add on regularity of the initial data in weighted Besov spaces, and obtain boundedness therein.
Due to the different structure of the nonlinear term $(\mathfrak{f},\mathfrak{g})$ compared to $(f,g)$, we need to improve the regularity of the low frequencies of $(a,u)$.
To that end, we will make use of a decay estimate from \cite{Danchin-Xu2017}.
For this result, we introduce the function space $X_p$ as the set of all pairs of functions $(a,u)$, where $a : [0,\infty) \times \R^d \to [0,\infty)$ is a scalar function and $u: [0,\infty) \times \R^d \to \R^d$ is a $d$-vector function, satisfying the following:
\begin{align*}
    & (a,u)^l \in \widetilde{C} (\R_{>0}; \db^{\frac{d}{2} - 1}_{2,1}) \cap L^1(\R_{>0}; \db^{\frac{d}{2} + 1}_{2,1}), 
    \\
    & a^h \in \widetilde{C} (\R_{>0}; \db^{\frac{d}{p} }_{p,1}) \cap L^1(\R_{>0}; \db^{\frac{d}{p} }_{p,1}),
    \\
    & u^h \in \widetilde{C}(\R_{>0}; \db^{\frac{d}{p}}_{p,1})
    \cap L^1(\R_{>0}; \db^{\frac{d}{p} + 1}_{p,1})
    .
\end{align*}
$X_p$ is then equipped with the obvious norm corresponding to the strong topologies for the above spaces.
The space $X_p$ is the original `critical space' used for the global existence theorems in \cites{danchin2016, Danchin-Xu2017}.

\begin{prop} {\rm (}\cite{Danchin-Xu2017}{\rm )}\label{dx2017}
    Let $d\geq2$ and $p \in [2, \text{min} \{4, 2d/(d-2)\} ]$, with $p\neq 4$ in the $d=2$ case.
    Assume without loss of generality that $P'(1) = \nu = 1$.
    Then there exists a constant $c = c(p,d, \mu, P)>0$ such that if
    \begin{align*}
        X_{p,0} \coloneqq 
        \norm{ (a_0,u_0) }^l_{ \db^{ \frac{d}{2}- 1}_{2,1} }
        +
        \norm{a_0}^h_{ \db^{\fdp }_{p,1} }
        +
        \norm{u_0}^h_{ \db^{ \fdp  - 1 }_{p,1} }
        \leq 
        c,
    \end{align*}
    then \eqref{31} has a unique global-in-time solution $(a,u)$ in $X_p$. 
    Furthermore, there exists a constant $C=C(p,d,\mu, P)>0$ such that
    \begin{align*}
        \norm{(a,u)}_{X_p} \leq C X_{p,0}.
    \end{align*}
    Also, there exists a constant $c_1$ such that if, in addition,
    \begin{align*}
        \norm{ (a_0,u_0) }^l_{ \db^{ -s_0 }_{2,\infty} } 
        \leq 
        c_1,
        \text{ where }
        s_0 \coloneqq d\biggc{ \frac{2}{p} - \frac{1}{2} },
    \end{align*}
    then we have a constant $C_1$ such that for all $t\geq 0$,
    \begin{align} \label{est of D}
        D_{p,\epsilon}(t) \leq C_1 \biggc{
        \norm{ (a_0,u_0) }^l_{ \db^{ -s_0 }_{2,\infty} } 
        +
        \norm{ (\N a_0,u_0) }^h_{ \db^{\fdp  - 1 }_{p,1} } 
        },
    \end{align}
    where the norm $D(t)$ is defined by 
    \begin{align*}
        D(t)
        \coloneqq
        \sup_{ s\in [\epsilon -s_0, d/2 + 1] }
        \norm{ 
        \langle 
        \tau
        \rangle^{
        %\frac{ s_0 + s }{ 2 } 
        (s_0+s)/2
        } 
        (a,u)
        }^l_{ L^\infty_t \db^s_{2,1} }
        +
        \norm{
        \langle 
        \tau
        \rangle^{ \fdp  + 1/2 - \epsilon}
        (\N a, u)
        }^h_{ \widetilde{L}^\infty_t \db^{ \fdp  - 1 }_{p,1}  }
        +
        \norm{ \tau \N u }^h_{ \widetilde{L}^\infty_t \db^{ \fdp  }_{p,1} },
    \end{align*}
    with $\epsilon >0$ taken sufficiently small.
\end{prop}

In particular, we will make use of the estimate \eqref{est of D} on the low-frequency part of $D_{p,\epsilon}(t)$ for $s = d/2$ 
%in the $p=2$ ($s_0 = d/2$) case 
to say that for all $t\geq0$,
\begin{align*}
    \norm{ (a,u)(t) }^l_{ \db^{\frac{d}{2}}_{2,1} }
    \leq 
    C  \langle t \rangle^{ - \fdp }  
    \biggc{
        \norm{ (a_0,u_0) }^l_{ \db^{ -s_0 }_{2,\infty} } 
        +
        \norm{ (\N a_0,u_0) }^h_{ \db^{ \fdp - 1 }_{p,1} } 
        }.
\end{align*}
We then obtain that, for $d\geq 3$ and $p \in [2,d)$, there exists a constant $C=C(p,d,\mu,P) >0$ such that, for all $t>0$,
\begin{align*}
    \norm{ (a,u) }^l_{ L^1_t \db^{\frac{d}{2}}_{2,1} } \leq C
    \biggc{
        \norm{ (a_0,u_0) }^l_{ \db^{ -s_0 }_{2,\infty} } 
        +
        \norm{ (\N a_0,u_0) }^h_{ \db^{ \fdp - 1 }_{p,1} } 
    }
    .
\end{align*}
The decay is just barely too slow to recover the $d=2$ case.

Finally, we prove the last result necessary for our main theorem. 
Essentially we add regularity of the initial data in weighted Besov spaces on top of the regularity assumed in Proposition \ref{dx2017}, and then apply the same steps as the proof of Proposition \ref{gl ext Danchin extended} to $(A_k,U_k)$ for their appropriate regularity exponents.

\begin{prop} \label{glob ests weighted} 
%{\bl (Glob Ests! weighted)}
    Let $d\geq 3$.
    %$k \in \{ 1, 2, ..., d \}$. 
    %$p\in[2,\text{min}\{d,2d/(d-2)\})$. 
    Assume
    %without loss of generality that 
    $P'(1) = \nu = 1$.
    Finally, let $(a_0,u_0)$ satisfy the conditions of Proposition \ref{dx2017}, and let $(a,u)$ be the associated solution to \eqref{31}.
    Then there exist a frequency cut-off constant $j_0 \in \Z$ and a small constant $c = c(d,\mu,P) >0$ such that, if $(a_0,u_0)$ 
    %and $(A_{k0},U_{k0})$ 
    satisfy 
    \begin{align*}
    \notag
        % S_{p,0} & \coloneqq  
        % \norm{(a_0,u_0)}_{L^\frac{p}{2} }
        % + \norm{(a_0,u_0)}^l_{\dot{B}^{\frac{d}{2} - 1}_{2,1}} 
        % + \norm{a_0}^h_{\db^{\frac{d}{p} + 1}_{p,1}} 
        % + \norm{u_0}^h_{\db^{\frac{d}{p}}_{p,1}}
        % \\
        % & \quad 
        % + \norm{(x_k a_0, x_k u_0)}^l_{\db^{\frac{d}{2}}_{2,1}  }
        % + \norm{x_k a_0}^h_{\db^{\frac{d}{p}+1}_{p,1}}
        % + \norm{x_k u_0}^h_{\db^{\frac{d}{p}}_{p,1}  }
        % \leq c,
        S_{0} \coloneqq &  
        \norm{(a_0,u_0)}^l_{\dot{B}^{\frac{d}{2} - 1}_{2,1}} 
        + \norm{a_0}^h_{\db^{\frac{d}{2} + 1}_{2,1}} 
        + \norm{u_0}^h_{\db^{\frac{d}{2}}_{2,1}}
        \\
        & + 
        \sum_{k=1}^d
        \Bigc{
        \norm{ ( x_k a_0, x_k u_0 ) }^l_{\db^{\frac{d}{2}}_{2,1} }
        +
        \norm{x_k a_0}^h_{ \db^{\frac{d}{2} + 1}_{2,1}}
        + \norm{x_k u_0}^h_{\db^{\frac{d}{2}}_{2,1}}
        }
        +
        \norm{(a_0,u_0)}_{ \db^{-\frac{d}{2}}_{2,\infty} }
        \leq c,
    \end{align*}
    then $(a,u)$ is in the space $S$. Also, there exists a constant $C = C(d,\mu, P, j_0)>0$ such that
    \begin{align*}
        \norm{(a,u)}_{S} \leq C S_{0}.
    \end{align*}
\end{prop}

\begin{pf}
By Proposition \ref{gl ext Danchin extended} we have global existence of the solution $(a,u)$ in $Y_{2}$, which contains $S$.
It remains only to show that the weighted norms are bounded. 
We will show boundedness of the weighted norms by considering the functions $(A_k,U_k)$ as solutions of \eqref{78x}.
We want
\begin{align*}
    \notag 
    S(t) \coloneqq & Y_2(t) +
    \sum_{k=1}^d
    \Bigc{
    \norm{ (A_k, U_k) }^l_{ L^\infty_t  \db^{ \frac{d}{2} }_{2,1}  }
    +
    \norm{U_k}^h_{L^\infty_t \db^{ \frac{d}{2} }_{2,1} }
    +
    \norm{A_k}^h_{ L^\infty_t \db^{ \fdtwo +1}_{2,1} }
    \\
    & +
    \norm{ (A_k, U_k) }^l_{ L^1_t  \db^{ \frac{d}{2} + 2}_{2,1}  }
    +
    \norm{U_k}^h_{L^1_t \db^{ \fdtwo + 2}_{2,1} }
    +
    \norm{A_k}^h_{ L^1_t \db^{ \fdtwo  + 1}_{2,1} }
    }
    \leq C S_{0},
\end{align*}
for all $t>0$.
We once again start by applying Lemma \ref{danchin prop 7} to the low-frequency norms, but this time we take $s=d/2$. We obtain
\begin{align} \label{weighted gl lin est for low freqs}
    \norm{ (A_k,U_k) }^l_{\widetilde{L}^\infty_t \db^{ \frac{d}{2} }_{2,1}}
    +
    \norm{ (A_k,U_k) }^l_{ L^1_t \db^{ \frac{d}{2} + 2}_{2,1}}
    \leq C
    \biggc{
    \norm{ (A_{k,0},U_{k,0} ) }^l_{ \db^{ \frac{d}{2}}_{2,1}}
    +
    \norm{ ( \mathfrak{f}, \mathfrak{g} ) }^l_{ L^1_t \db^{  \frac{d}{2} }_{2,1}}
    }.
\end{align}

Next, observing from \eqref{78x} that $\P U_k$ solves a heat equation, we apply Proposition \ref{thm1} to get
\begin{align} \label{weighted gl lin est for Pu }
    \norm{ \P U_k }^h_{\widetilde{L}^\infty_t \db^{\frac{d}{2}}_{2,1} }
    +
    \mu \norm{ \P U_k }^h_{ L^1_t \db^{\frac{d}{2} + 2}_{2,1} }
    \leq
    C \biggc{
    \norm{ \P U_{k,0} }^h_{  \db^{\frac{d}{2}}_{2,1} }
    +
    \norm{ \P \mathfrak{g} }^h_{ L^1_t \db^{\frac{d}{2} }_{2,1} }
    }.
\end{align}
Mimicking Proposition \ref{gl ext Danchin extended}, we define the `weighted effective velocity'
\begin{align*}
    W_k \coloneqq \N (-\Delta)^{-1} (A_k -\div(U_k) ),
\end{align*}
which solves 
\begin{align*}
    \pt_t W_k - \Delta W_k 
    =
    \N(-\Delta)^{-1} (\mathfrak{f} -\div(\mathfrak{g}) )
    + W_k - (-\Delta)^{-1} A_k,
\end{align*}
and so we may apply Proposition \ref{thm1} to $W_k$ as well to get
\begin{align} \label{weighted gl lin est for w}
    &  
    \norm{ W_k }^h_{ \widetilde{E}^{ \fdtwo }_2 (t) } 
    % \coloneqq
    % \norm{ W_k }^h_{\widetilde{L}^\infty_t \db^{d/p}_{p,1} }
    % +
    % \norm{ W_k }^h_{ L^1_t \db^{d/p + 2}_{p,1} }
    % \\
    % & 
    \leq C \biggc{
    \norm{W_{k,0}}^h_{\db^{  \fdtwo }_{2,1}}
    +
    \norm{ \mathfrak{f} - \div(\mathfrak{g}) }^h_{ L^1_t \db^{\frac{d}{2} - 1 }_{2,1}} 
    +
    \norm{ W_k - (-\Delta)^{-1} \N A_k }^h_{ L^1_t \db^{\frac{d}{2} }_{2,1}} 
    }
\end{align}
Next, noting that 
\begin{align*}
    \pt_t A_k + \div(A_k u) + A_k
    =
    -\div(W_k) + au_k + u_k,
\end{align*}
we may bound $A_k$, similarly to the proof of Proposition \ref{thm2}, to get
\begin{align} \label{weighted gl lin est for a }
    & \notag 
    \norm{ A_k }^{h}_{\widetilde{L}^\infty_t \db^{ \fdtwo +1}_{2,1}}
    +
    \norm{ A_k }^{h}_{ L^1_t \db^{ \fdtwo +1}_{2,1} }
    \\
    & \leq C \notag
    \Bigc{
    \norm{ A_{k,0} }^h_{ \db^{  \fdtwo + 1 }_{2,1} }
    +
    \norm{ \div(W_k) }^h_{ L^1_t \db^{ \fdtwo +1}_{2,1} }
    +
    \norm{ au_k + u_k }^h_{ L^1_t \db^{ \fdtwo +1}_{2,1} }
    \\
    & \quad  \quad  \quad  \quad +
    \int^{t}_{0}
    (
    \norm{\N u}_{  \db^{  \fdtwo   }_{2,1} }
    \norm{A_k}_{  \db^{  \fdtwo  + 1 }_{2,1} }
    +
    \norm{\N u}_{  \db^{  \fdtwo  + 1 }_{2,1} }
    \norm{A_k}_{ \db^{  \fdtwo  }_{2,1} }
    )
    \dd \tau
    }.
\end{align}
Then, similarly to Proposition \ref{gl ext Danchin extended}, we take a sufficiently large cut-off constant $j_0$ to combine \eqref{weighted gl lin est for low freqs}, \eqref{weighted gl lin est for Pu }, \eqref{weighted gl lin est for w}, and \eqref{weighted gl lin est for a } and obtain
\begin{align*}
    S_p(t) \leq \notag
    &
    C 
    \Biggc{
    S_{p,0}
    +
    \int^t_0
    (
    \norm{ ( \mathfrak{f}, \mathfrak{g} ) }^l_{ \db^{ \frac{d}{2} }_{2,1} }
    +
    \norm{ \mathfrak{f} }^h_{ \db^{  \fdtwo  - 1 }_{2,1} }
    +
    \norm{ \mathfrak{g} }^h_{ \db^{  \fdtwo  }_{2,1} }
    +
    \norm{au}^h_{ \db^{\frac{d}{2} + 1}_{2,1} }
    +
    \norm{u}^h_{ \db^{\frac{d}{2} + 1}_{2,1} }
    \\ &
    \quad \quad \quad \quad 
    +
    \norm{\N u}_{  \db^{  \fdtwo   }_{2,1} }
    \norm{A_k}_{ \db^{  \fdtwo  + 1 }_{2,1} }
    +
    \norm{\N u}_{  \db^{  \fdtwo + 1 }_{2,1} }
    \norm{A_k}_{  \db^{  \fdtwo }_{2,1} }
    )
    \dd \tau
    }.
\end{align*}
We consider briefly the low frequencies of the nonlinear terms. Expanding $\mathfrak{f}$,
\begin{align*}
\norm{ \mathfrak{f} }^l_{L^1_t \db^{\frac{d}{2}}_{2,1} } 
=
\norm{ 
-A_k \div(u) - u \cdot \N A_k + au_k + u_k
}^l_{L^1_t \db^{\frac{d}{2}}_{2,1} } .
\end{align*}
For the `lone term' $u_k$, we simply apply the low-frequency estimate in Proposition \ref{dx2017}.
For the remaining terms, we may apply the same product estimates as for $f$ in \cite{danchin2016}, since the regularity exponent is already $d/2$ (one higher than it was for $f$).
We may similarly treat $\mathfrak{g}$.
% {\bl (It might be best to check this part of the proof together, or for me to write it out in more detail in this paper, although that could take up a lot of space. How important is the 20-page limit?)}

We thus obtain 
\begin{align*}
    S(t) \leq C \biggc{
    S_{0} + S(t)^2
    },
\end{align*}
and again thanks to our local existence result Proposition \ref{local exist thm hi reg} and smallness of $S_{0}$, this allows us to extend to global existence of $(a,u)$ in $S.$
\end{pf}

\bigskip

\noindent {\bf Acknowledgements.} The author would like to thank their advisor, Professor Tsukasa Iwabuchi of Tohoku University, for their kind support and assistance.

\phantom{}

\noindent {\bf Conflict of interest statement.} On behalf of all authors, the corresponding author states that there is no conflict of interest. 

\phantom{}

\noindent {\bf Data availability statement.} This manuscript has no associated data. 

\phantom{}

\bibliographystyle{unsrt} \bibliography{sample}

\end{document}